\documentclass[a4paper,11pt,leqno]{article}
\pdfoutput=1
\usepackage[T2A,T1]{fontenc}
\usepackage[utf8]{inputenc}
\usepackage{CJKutf8} 
\usepackage[russian,greek,english]{babel}
\usepackage{amsmath,amsfonts,amsthm,amssymb}
\usepackage[left=3cm, right=3cm, top=3cm, bottom=3.5cm]{geometry}
\usepackage{graphicx}
\usepackage{enumerate}

\usepackage[hyphens]{url}
\usepackage{hyperref}

\hypersetup{
 pdftitle={S1 Porto},
 pdfauthor={Bernold Fiedler},
 colorlinks=true,
 linkcolor=blue,
 citecolor=blue,
 filecolor=blue,
 urlcolor=blue}

 \usepackage{afterpage}

\usepackage{wrapfig} 
\usepackage[format=plain,labelfont=bf,font=small]{caption}
\usepackage{xcolor}
\usepackage[arrow, matrix, curve]{xy}
\usepackage{float}

\usepackage{caption}
\usepackage{tabulary}
\usepackage{array}
\newcolumntype{N}[1]{>{\centering\arraybackslash}m{#1}}

\newcommand{\transv}{\mathrel{\text{\tpitchfork}}}
\makeatletter
\newcommand{\tpitchfork}{%
  \vbox{
    \baselineskip\z@skip
    \lineskip-.52ex
    \lineskiplimit\maxdimen
    \m@th
    \ialign{##\crcr\hidewidth\smash{$-$}\hidewidth\crcr$\pitchfork$\crcr}
  }%
}
\makeatother
\usepackage{latexsym}

\usepackage[notref,notcite,color, final 
]{showkeys}

\definecolor{refkey}{rgb}{1,0,0}
\definecolor{labelkey}{rgb}{1,0,0}

  \mathchardef\ordinarycolon\mathcode`\:
  \mathcode`\:=\string"8000
  \begingroup \catcode`\:=\active
    \gdef:{\mathrel{\mathop\ordinarycolon}}
  \endgroup
  
\theoremstyle{plain}
\newtheorem{theorem}{Theorem}[section]
\newtheorem{lemma}[theorem]{Lemma}

\newtheorem{corollary}[theorem]{Corollary}

\newtheorem{definition}[theorem]{Definition}

\newcommand\eps{\varepsilon}
\newcommand\mi{\mathrm{i}}

\renewcommand\rho{\varrho}
\renewcommand\phi{\varphi}
\renewcommand\Re{\mathrm{Re}\,}

\newcommand{\N}{\mathbb{N}}
\newcommand{\Z}{\mathbb{Z}}
\newcommand{\R}{\mathbb{R}}

\def\cA{\mathcal{A}}
\def\cE{\mathcal{E}}
\def\cR{\mathcal{R}}
\def\cP{\mathcal{P}}
\def\cF{\mathcal{F}}
\def\cH{\mathcal{H}}
\def\cC{\mathcal{C}}
\def\cD{\mathcal{D}}
\def\cS{\mathcal{S}}
\def\cT{\mathcal{T}}

\def\cM{\mathcal{M}}
\def\cN{\mathcal{N}}

\newcommand{\sss}[1]{\subsection{#1}}

\def\cAP{\mathcal{A}_f^\mathcal{P}}
\def\AT{\cA_\cT^\cD}
\def\cEP{\mathcal{E}_f}
\def\cFP{\mathcal{F}_f^\mathcal{P}}
\def\cRP{\mathcal{R}_f^\mathcal{P}}
\def\cHP{\mathcal{H}_f^\mathcal{P}}
\def\cCP{\mathcal{C}_f^\mathcal{P}}
\def\CT{\cC_\cT^\cD}

\def\cgAP{\mathcal{A}_g^\mathcal{P}}

\def\cgCP{\mathcal{C}_g^\mathcal{P}}

\def\cAN{\mathcal{A}_g^\mathcal{N}}
\def\cEN{\mathcal{E}_g}
\def\cFN{\mathcal{F}_g^\mathcal{N}}
\def\cHN{\mathcal{H}_g^\mathcal{N}}
\def\cCN{\mathcal{C}_g^\mathcal{N}}

\newcommand{\bS}{\mathbb{S}}
\newcommand{\bO}{\mathbb{O}(2)}
\newcommand{\bSO}{\mathbb{SO}(2)}

\def\SSOP{\mathrm{Sturm}^\cP(u,u_x)}
\def\SOP{\mathrm{Sturm}^\cP(u,u_x^2)}
\def\SHP{\mathrm{Sturm}^\cP(u)}
\def\SSON{\mathrm{Sturm}^\cN(u,u_x)}
\def\SON{\mathrm{Sturm}^\cN(u,u_x^2)}
\def\SHN{\mathrm{Sturm}^\cN(u)}
\def\SxN{\mathrm{Sturm}^\cN(x,u,u_x)}
\def\SO{\mathrm{Sturm}(u,u_x^2)}
\def\SH{\mathrm{Sturm}(u)}

\def\id{\mathop{\mathrm{id}}\nolimits}

\def\Re{\mathop{\mathrm{Re}}\nolimits}
\def\sign{\mathop{\mathrm{sign}}\nolimits}
\def\codim{\mathop{\mathrm{codim}}\nolimits}
\def\dim{\mathop{\mathrm{dim}}\nolimits}

\begin{document}

\title{\LARGE{\textbf{Classifications of global attractors \\[2mm]
for $\bS^1$-equivariant parabolic equations:\\
a survey\\}}}

\author{
\\[1em]
\Large{\emph{Dedicated, with due admiration, to}}\\[2mm]
\Large{\emph{Mike Field, Marty Golubitsky, John Guckenheimer, and Ian Stewart, }}\\[2mm]
\Large{\emph{on the occasion of their 80th birthdays}}
\\ 
\\[1em]
Carlos Rocha*, Bernold Fiedler**, Alejandro L\'{o}pez-Nieto***
\\[2cm]
}

\date{\small{version of \today}}
\maketitle
\thispagestyle{empty}
\vfill
\setlength{\parindent}{0cm}
*\\
        Centro de An\'alise Matem\'atica, Geometria e Sistemas Din\^amicos,\\
				Instituto Superior T\'ecnico, Universidade de Lisboa\\
        Av. Rovisco Pais 1, 1049--001 Lisbon, PORTUGAL\\
        \url{crocha@tecnico.ulisboa.pt}\\
        \url{http://camgsd.tecnico.ulisboa.pt}\\[5mm]
**\\
        Institut f\"ur Mathematik\\
        Freie Universit\"at Berlin\\
        Arnimallee 7, D--14195 Berlin, GERMANY\\
        \url{fiedler@math.fu-berlin.de}\\
        \url{http://dynamics.mi.fu-berlin.de}\\[5mm]
***\\
	Department of Mathematics\\
	National Taiwan University\\
	No. 1, Sec. 4, Roosevelt Road,\\
	10617 Taipei, TAIWAN \\
        \url{alopez@ntu.edu.tw}\\
        \url{https://sites.google.com/view/alnalejandro}\\


\newpage
\pagestyle{plain}
\pagenumbering{roman}
\setcounter{page}{1}

\begin{abstract}
\noindent
We survey the global dynamics of semiflows generated by scalar semilinear parabolic equations which are $\bSO$ equivariant under spatial shifts of $x\in \bS^1=\mathbb{R}/2\pi\mathbb{Z}$, i.e.
\begin{equation}
\label{aPDEP}
u_t = u_{xx} + f(u,u_x),\qquad x\in \bS^1.
\end{equation}
For dissipative $C^2$ nonlinearities $f$, the semiflow \eqref{aPDEP} possesses a compact global attractor $\cA=\cAP$ which we call Sturm attractor. 
The Sturm attractor $\cAP$ decomposes as
\begin{equation*}
\label{acAP}
\cAP=\cEP\cup\cFP\cup\cRP\cup\cHP,
\end{equation*}
where $\cHP$ denotes heteroclinic orbits between distinct elements of spatially homogeneous equilibria $\cEP$, rigidly rotating waves $\cRP$ and, as their non-rotating counterparts, frozen waves $\cFP$.
We therefore represent $\cAP$ by its \emph{connection graph} $\cCP$, with vertices in $\cEP\,,\cFP,\cRP$ and edges $\cHP$. 
Under mild hyperbolicity assumptions, the directed graphs $\cCP$ are finite and transitive.
For illustration, we enumerate all 21 connection graphs $\cCP$ with up to seven vertices.
The result uses a lap signature of period maps associated to integrable versions of the steady state ODE of \eqref{aPDEP}. 

\noindent
Our results are based on a comparison with associated $C^2$ dissipative Neumann PDEs on the half-interval $x\in(0,\pi)$.
After a suitable homotopy from $f$ to spatially reversible nonlinearities $g(u,-p)=g(u,p)$, the Neumann problem gains a \emph{hidden symmetry} $\bO$, spatially.
Moreover, the Neumann dynamics become gradient-like.
This leads to isomorphic connection graphs
\begin{equation*}
\label{acCPN}
\cCP\cong\cCN/\sim\,.
\end{equation*}
Here $\sim$ collapses the two shifted Neumann copies of each frozen or rotating wave to a single vertex.
As an example, we freeze and reconstruct the connection graph of the Vas tulip attractor, known from delay differential equations, in setting \eqref{aPDEP}.
\end{abstract}

\vspace{1cm}
{\small{
\tableofcontents
\listoffigures
\listoftables
}}


\newpage
\pagenumbering{arabic}
\setcounter{section}{-1}
\setcounter{page}{1}
\setcounter{equation}{0}
\vspace{1cm}

\section{Nontechnical summary}

\textbf{Global attractors, as introduced by Ladyzhenskaya \cite{lady} in 1972 for the two-dimensional Navier-Stokes equations of fluid dynamics, are a concept to capture the long-time dynamics of infinite-dimensional dynamical systems.
Morse-Smale systems, as introduced by Palis \cite{Palis} in 1969, are marked by local structural stability under small perturbations.
Our present survey presents a global classification of all Morse-Smale global attractors for partial differential equations (PDEs) of translation-invariant reaction-drift-diffusion type, in one space dimension and under spatially periodic boundary conditions.
This revisits mathematically rigorous work of the last 40 years.
In particular, we present and discuss
\begin{enumerate}[(i)]
    \item a symbolic encoding of the global Morse-Smale attractors;
    \item the classification of their global heteroclinic connection graphs among rotating waves, frozen waves, and homogeneous equilibria;
    \item the role of hidden symmetries under passage to non-periodic Neumann boundary conditions, i.e., to a gradient-like Morse case which freezes any time-periodic solutions, 
\end{enumerate}
For illustration, we enumerate all 21 examples with up to a total of 7 homogeneous equilibria and rotating/frozen waves.
We conclude with comments on some mysteriously analogous results for certain differential equations with a time-lagged argument, including the celebrated Vas tulip.}

\section{Introduction}\label{Intro}

\numberwithin{equation}{section}
\numberwithin{figure}{section}
\numberwithin{table}{section}

Dynamical systems provide an essential tool for the analysis of applications, e.g. in engineering, physics, chemistry, biology, sociology, etc. 
Infinite-dimensional applied models include semiflows generated by partial differential or delay differential equations.
We survey some results from the last four decades on scalar parabolic PDEs of the form
\begin{equation} \label{101} 
u_t = u_{xx} + f(u,u_x) \,, \quad x\in \bS^1=\R/2\pi\Z \,.
\end{equation}
We assume the nonlinearity $f{:}\ \R^2\rightarrow \R$ is $C^2$-smooth and \emph{dissipative}. 
By smoothness, \eqref{101} generates a local semiflow $\cS(t){:}\ X\rightarrow X, t\ge 0$, in the Sobolev space 
$X=H^s(\bS^1),\, s>3/2$.

Note $\bSO\cong(\bS^1,+)$ \emph{equivariance} of \eqref{101} under $x$-shifts of $u$ by any fixed $\xi\in\bS^1$.
Indeed, $u(t,x+\xi)$ is a solution whenever $u(t,x)$ is.
For celebrated background on differentiable dynamics in general, equivariant dynamics in particular, and applications, see the monumental works by our jubilees Mike Field, Marty Golubitsky, John Guckenheimer, and Ian Stewart, like \cite{field89, field15, golste88, golste23, gu06, guho83}.

\sss{Global attractors}\label{A}

Dissipativity requires that any solution $u(t,\cdot)=\cS(t)u_0$ enters some fixed large ball in $X$, for all sufficiently large times $t\geq t_0(u_0)$.
In particular, the semiflow $\cS(t)$ is global, for all $t\geq 0$.
Sufficient conditions for dissipativity, which we impose throughout, are
\begin{equation}
\label{diss}
f(v,0)v<0\quad \textrm{and} \quad |f(v,p)|\leq C(|v|)(1+p^2)\,,
\end{equation}
for all large $|v|$ and some continuous function $C(|v|)$; see \cite{mana97}.
In view of homogeneous equilibria, the sign condition is also necessary.

For dissipative $f$, the PDE \eqref{101} possesses a nonempty compact \emph{global attractor} $\cAP\subset X$. 
The superscript $\cP$ indicates periodic boundary conditions. 
The global attractor $\cAP$ attracts all bounded sets. 
It consists of all solutions $u(t,\cdot)$ of \eqref{101} which exist and remain uniformly bounded for all $t\in\mathbb{R}$, forward \emph{and} backward.
See \cite{hen81, paz83} for some background on semiflows, and \cite{bavi92, hal88, hamaol02, lady} for global attractors. 
For \eqref{101} more specifically, see  \cite{anfi88, fm-p, mat88, mana97}, and \cite{firowo04,firowo12b, roc25} for details. 

By definition, the global attractor $\cAP$ contains all spatially \emph{homogeneous equilibria} $u(t,x)\equiv e\in\cE,\ f(e,0)=0$. 
Dissipativity implies that there exists a smallest zero $e=\underline{e}\in\R$ and a largest zero $\overline{e}$, possibly identical.
We assume nondegeneracy $f_v(e,0)\neq0$ for all $e\in\cE$.
By dissipativity of $f$, the number $n$ of homogeneous equilibria $e_j$ is then odd.

Due to $\bS^1$-equivariance, spatially nonhomogeneous \emph{rotating waves} $u(t,x) =$ $v(x-ct)\in\cR$, which rotate at constant wave speeds $c\neq0$, are another option. 
Variedly, such translating PDE solutions are also called \emph{traveling waves} and, particularly in integrable context, \emph{solitary waves}.
Indeed, these are relative equilibria: time action is equivalent to spatial group action of $x$-shift by $\xi=-ct$.
Specifically, rotating waves correspond to solutions $v=v(x)$ of the co-rotating ODE
\begin{equation} 
\label{103} 
0 = v_{xx} + f(v,v_x) +cv_x\,, \quad x\in \bS^1 \,.
\end{equation}
In other words, the nonhomogeneous wave profile $v$ possesses \emph{minimal spatial period} $T=2\pi/\ell$, for some integer $\ell>0$, and \emph{minimal time-period} $\Theta=T/c$.

For vanishing speed parameter $c=0$, we obtain nonhomogeneous equilibria $u(t,x)=v(x)$, which solve
\begin{equation} \label{102} 
0 = v_{xx} + f(v,v_x) \,, \quad x\in \bS^1 \,.
\end{equation}
As usual, we rewrite the second order ODE \eqref{102} as a system
\begin{equation}
\label{vpf}
\begin{aligned}
v_x&=p\,,\\
p_x&=-f(v,p)\,.
\end{aligned}
\end{equation}
Substituting $f$ by $f+cp$, we subsume \eqref{103} under the same ODE form.
The homogeneous PDE equilibria $e\in\cE$ remain ODE equilibria $(e,0)$, equivalently.

Next consider any nonstationary periodic ODE orbit $(v,p)$ with minimal period $T=T_f(v)>0$.
Extrema are located on the horizontal $v$-axis $p=0$.
In particular, they are nondegenerate.
By the Jordan curve theorem, there are only two extremal values, for each periodic orbit $v=v(x)$:
\begin{equation}
\label{minmax}
\underline{v}:=\min v \quad\textrm{and}\quad \overline{v}:=\max v\,.
\end{equation}
Therefore $T=T_++T_-$\,, where $T_+$ measures the $x$-interval spent above the horizontal $v$-axis, with $p=p_+(v)>0$, and $T_-$ accounts for the lower counterpart $p=p_-(v)<0$.
Explicitly,
\begin{equation}
\label{Tpm}
T_\pm=\pm\int_{\underline{v}}^{\overline{v}} \,1/p_\pm(v)\,dv\,.
\end{equation}

The nonhomogeneous PDE equilibria $v$, specifically, become nonstationary periodic ODE orbits $(v,p)$ with minimal  period
\begin{equation}
\label{TlP}
0<T_f(v)=2\pi/\ell\,,
\end{equation}
for some positive integer $\ell\in\N$.
The $\bS^1$-action by spatial shifts $\xi$ ensures that each such $v$ generates a circle's worth of PDE equilibria.
We call such circles $\bS^1v$ \emph{frozen waves} $\cF$, to emphasize that they occur \emph{as a single symmetry related group orbit,} just as the time orbits of rotating waves do; see \cite{fiedlerhabil}.

Keeping track of $f$ and the periodic boundary conditions $\cP$, the global attractor $\cAP$ in fact decomposes as
\begin{equation}
\label{cAP}
\cAP=\cEP\cup\cFP\cup\cRP\cup\cHP\,.
\end{equation}
Here $\cHP$ denote \emph{heteroclinic}, alias \emph{connecting}, time orbits $u(t,x)$ between distinct \emph{sources} $v_1=\lim_{t\rightarrow -\infty} u(t,\cdot)$ and \emph{targets} $v_2=\lim_{t\rightarrow +\infty} u(t,\cdot)$ in $\cEP\cup\cFP\cup\cRP$\,.
In symbols,
\begin{equation}
\label{leadsto}
u(t,\cdot):\quad v_1\leadsto v_2\,.
\end{equation}
This follows from $\bS^1$-equivariance \cite{mana97} and a Poincaré-Bendixson theorem \cite{fm-p} for \eqref{101}.

\sss{Hyperbolicity}\label{Hyp}

Hyperbolicity of sources and targets $v\in \cEP\cup\cFP\cup\cRP$ is therefore a central notion to the global analysis of the global attractor decomposition \eqref{cAP}.
For equilibria $v\in\cEP\cup\cFP$, we have to consider the eigenvalues $\lambda$ of the linearization
\begin{equation} \label{110}
\lambda w = w_{xx} + f_p(v,v_x) w_x + f_u(v,v_x) w \,, \quad x\in \bS^1 \,.
\end{equation} 
At homogeneous equilibria $v\equiv e\in\cEP$, we obtain trigonometric eigenfunctions with explicit spectrum, alias dispersion relation, $\lambda=\lambda_k=-k^2\pm\mi k f_p(e,0)+f_u(e,0),\ k\in\Z$.
Since homogeneous equilibria $e$ are invariant under the $\bS^1$-action on $X$, this Fourier decomposition into eigenspaces $E_k=\langle \cos kx, \sin kx \rangle$ coincides with the decomposition of $X$ into $\bS^1$-irreducible subspaces. 
Note $\dim E_0=1$, and $\dim E_k=2$ for integer $k>0$.
\emph{Hyperbolicity} of $e$ requires 
\begin{equation}
\label{hyp}
\Re \lambda_k\neq 0\,,
\end{equation} 
for all $k$. 
The \emph{unstable dimension}, alias the \emph{Morse index} $i(e)$, then counts the eigenvalues with $\Re \lambda_k>0$, repeated with algebraic multiplicities.
Indeed, 
\begin{equation}
\label{ie}
i(e)=\dim W^u(e)=\codim W^s(e) 
\end{equation}
is the dimension of the \emph{unstable manifold} $W^u(e)$ of $e$, and the codimension of the \emph{stable manifold} $W^s(e)$.
For $f_u(e,0)<0$, note that $e$ is ``hyperbolically'' stable: $i(e)=0$.
For $k^2<f_u(e,0)<(k+1)^2,\ k\geq0$, we obtain hyperbolicity with odd unstable dimension $i(e)=2k+1$.

Successive homogeneous equilibria $e_j\in\cEN$, for example, are ODE saddles $(e_j,0)$, for $f_v(e_j,0)<0$, and ODE sinks, sources, or centers, for $f_v(e_j,0)>0$, alternatingly in $j$.
ODE saddles $e$ are PDE-stable, and ODE sinks/sources/centers $e$ are PDE-unstable.
The minimal and maximal equilibria $\underline{e}=e_1$ and $\overline{e}=e_n$ are PDE-stable, by dissipativity, and hence ODE saddles.
In particular, the number $n=|\cE_f|$ of homogeneous equilibria $e_j$ is odd.

All frozen or rotating waves $v$ which we discuss next, in contrast, will be unstable by Sturm-Liouville comparison of \eqref{110} with the nonconstant eigenfunction $w=v_x$ for the trivial eigenvalue $\lambda=0$.
The same argument extends to Neumann boundary conditions discussed below; see \eqref{PDEg}.

For \emph{hyperbolicity of frozen waves} $v=v(x)\in\cFP$, the trivial eigenfunction $w=v_x \not\equiv 0$ poses an obstacle.
Indeed, $v_x$ is the tangent to the group orbit $\bS^1v$, which constitutes any frozen wave.
Therefore, we define hyperbolicity of frozen waves by \eqref{hyp}, except for the trivial eigenvalue $\lambda=0$ which we require to be algebraically simple.
This amounts to \emph{normal hyperbolicity} of the frozen wave $\bS^1v$.
The \emph{Morse index} $i(v)$ again counts eigenvalues $\Re \lambda >0$ with their algebraic multiplicities.
In particular, $i(v)=\dim W^{uu}(v)$ provides the dimension of the \emph{strong unstable manifold} $W^{uu}(v)$ of $v$.
Similarly, $\codim W^{ss}(v)=i(v)+1$ is the codimension of the \emph{strong stable manifold} $W^{ss}(v)$.
The \emph{unstable} and \emph{stable manifolds} $W^u$ and $W^s$ of the frozen wave group orbit $\bS^1v$ are defined as the group orbits $\bS^1W^{uu}(v)$ and $\bS^1W^{ss}(v)$, respectively.
This view point agrees with our decision to consider group orbits $\bS^1v$ as frozen waves, rather than individual equilibria $v$ on them.
For the group orbits,
\begin{equation}
\label{iv}
\begin{aligned}
\dim W^u &= i(v)+1\,;\\
\codim W^s &= i(v)\,.
\end{aligned}
\end{equation}

For \emph{hyperbolicity of rotating waves} $u=v(x-ct)\in\cRP$, the above distinctions are commonplace.
The role of eigenvalues $\lambda$ is then played by Floquet exponents, with associated Floquet multipliers $\mu=\exp(\lambda \Theta)$.
As before, $\Theta=T/c$ is the minimal time-period of $v$ with minimal spatial period $T$.
When counting multiplicities of Floquet exponents, we therefore have to identify Floquet exponents $\lambda$ which only differ by integer multiples of $2\pi\mi/\Theta$ and therefore give rise to the same Floquet multiplier $\mu$.
Hyperbolicity of time-periodic orbits requires \eqref{hyp}, except for an algebraically simple trivial Floquet exponent $\lambda=0$.
The analogy to frozen waves becomes particularly evident if we pass to co-rotating dependent variables $\tilde{u}(t,x):=u(t, x+ct)$.
Then \eqref{101} becomes
\begin{equation}
\label{101c}
\tilde{u}_t=\tilde{u}_{xx}+f(\tilde{u},\tilde{u}_x)+c\tilde{u}_x \,, \quad x\in \bS^1=\R/2\pi\Z \,.
\end{equation}
Note how the rotating wave $u$ becomes a frozen wave $\tilde{u}(t,x)=v(x)$, and compare with the equilibrium equation \eqref{103}.
Hyperbolicity and Morse indices do not change, because co-rotation only affects the imaginary parts of the spectrum $\lambda$.
In particular, the unstable manifold $W^u$ of the rotating wave $u=v(x-ct)$ coincides with the unstable manifold $W^u$ of the frozen wave $\tilde{u}=\bS^1 v$, and satisfies \eqref{iv}.
The strong unstable manifolds $W^{uu}$ of the frozen wave track convergence in phase towards the original rotating wave, for $t\searrow-\infty$.

In addition to $C^2$-regularity and dissipativity, we henceforth assume hyperbolicity of all homogeneous equilibria, frozen and rotating waves in \eqref{110}--\eqref{iv}, i.e.
\begin{equation}
\label{HYP}
\textbf{all elements of } \cEP\cup\cFP\cup\cRP \textbf{ are hyperbolic.}
\end{equation}

As a variant of PDE \eqref{101}, we will consider the associated \emph{Neumann problem} on the half-interval $x\in(0,\pi)$, i.e.
\begin{equation}
\label{PDEg} 
u_t = u_{xx} + g(u,u_x) \,, \quad u_x=0 \textrm{ at } x=0,\pi\,.
\end{equation}
See for example \cite{firo96} based on \cite{furo91}.
For some early results on Dirichlet boundary conditions which were technically limited to the Hamiltonian case $f=f(u)$, see \cite{brfi88, brfi89}.

In both, the periodic and the Neumann setting, the spatially reversible case
\begin{equation}
\label{grev}
g(v,-p)=g(v,p),\quad\textrm{for all } v,p\in\R
\end{equation}
deserves special attention.
Note how spatial reversibility \eqref{grev} makes the original PDE semiflow \eqref{101} for $f=g$ equivariant under $\bO$. 
The $\bO$-action is generated by the $x$-shifts of $\bSO$ and the \emph{flip}
\begin{equation}
\label{flip}
\kappa u(t,x):=u(t,2\pi-x)) 
\end{equation}
through $x=0,\pi$.
In group jargon, $\bO=\bSO\rtimes\langle\kappa\rangle$ is the semidirect product of the normal subgroup $\bSO$ with the group $\Z_2\cong\langle\kappa\rangle$ generated by the flip.

\sss{Sturm classes}\label{Sturm}

We collect our respective assumptions in the following definition.
\begin{definition}\label{defsturm}
We denote the $\bSO$ \emph{Sturm class} of $C^2$ nonlinearities $f$ in \eqref{101}, which satisfy dissipativity and hyperbolicity assumptions \eqref{diss} and \eqref{HYP}, by
\begin{equation}
\label{sturmP}
f\in\SSOP\,.
\end{equation}
If a nonlinearity $f\in\SSOP$ satisfies \eqref{grev}, in addition, we obtain the \emph{spatially reversible} $\bO$ \emph{class}
\begin{equation}
\label{sturmP2}
f\in\SOP\,.
\end{equation}
The \emph{Hamiltonian class} requires $f=f(u)\in\SOP$ to be independent of $p$ altogether:
\begin{equation}
\label{sturmP0}
f\in\SHP\,.
\end{equation}
The Hamiltonian class reduces the steady state ODEs \eqref{102} and \eqref{vpf} to the Hamiltonian case of a second order pendulum.

Analogous notation with superscripts $\cN$ will refer to the Neumann variant \eqref{PDEg} of PDE \eqref{101} on the half-interval $x\in(0,\pi)$.
In fact, we slightly extend the above terminology in the Neumann case, only.
Let
\begin{equation}
\label{sturmNx}
g\in\SxN
\end{equation}
denote the class of nonlinearities $g=g(x,u,u_x)$ in PDE \eqref{PDEg} which are also allowed to depend on $x$, but remain $C^2$, dissipative, and hyperbolic.
\end{definition}
For example, note the inclusions $\SHP\subset\SOP\subset\SSOP$, and similarly $\SHN\subset\SON\subset\SSON\subset\SxN$.
Moreover $\SOP=\SON$ and $\SHP=\SHN$.
Indeed, restriction of solutions $u$ of \eqref{101} to the time-invariant isotropy subspace $\{\kappa u=u\}$ shows ``$\subseteq$'' in the last two equalities.
In view of dissipativity conditions \eqref{diss} above and hyperbolicity result \eqref{hypT} below, which do not distinguish $\cP$ or $\cN$, equality follows.
It is therefore safe to omit superscripts $\cP,\cN$ in these two classes and write $\SH\subset\SO$, instead.

\sss{Zero numbers}\label{Zero}

The perhaps most powerful technical feature of PDE \eqref{101}, in addition to $\bS^1$-equivariance, is the dropping behavior of the even \emph{zero number} 
\begin{equation}
\label{z}
z^\cP:\quad X\setminus\{0\}\rightarrow \{0,2,4,\ldots\}\,.
\end{equation}
Specifically, $z^\cP(w)$ counts the number of strict sign changes of $x\mapsto w(x)\in X\setminus\{0\}\subset C^1$.
Let $w(t,x):=u_2(t,x)-u_1(t,x)$ denote the difference of any two nonidentical solutions $u_1,u_2$ of \eqref{101}. 
Then 
\begin{equation} 
\label{109}
t \mapsto z^\cP(w(t,\cdot)) \searrow
\end{equation} 
is nonincreasing with $t$.
More precisely, $z^\cP$ becomes finite, as soon as $t>0$, and drops strictly at any $t_0>0$ for which $x\mapsto w(t_0,x)$ possesses a multiple zero $w=w_x=0$ at some $x=x_0\in\bS^1$.
The same properties extend to Neumann or other separated boundary conditions.
See \cite{ang88, mat82}.

With the original idea going back to Charles-François Sturm \cite{sturm} in 1836, zero number dropping \eqref{109} is based on a detailed analysis of linear equations
\begin{equation} \label{wlin}
w_t = w_{xx} + b(t,x) w_x + a(t,x)w \,, \quad x\in \bS^1 \,, \ t>0 \,, \quad w(0,x) = w_0(x) \,,
\end{equation}
with coefficients $a,b$ depending on $u_1,u_2$.
This is the reason why we attach the name of Sturm to the above classes of nonlinearities, and to the resulting global attractors.

Linearization of the semiflow $(\cS(t)u_0)(x):=u(t,x)$ generated by \eqref{101} along the 
solution $u(t,x)$ leads to an infinitesimal version of \eqref{wlin} with $b= f_p(u,u_x),\  a=f_u(u,u_x)$.
Linearization provides a linear nonautonomous evolution $\cT(t,\tau),\ t\ge\tau$, which elucidates the asymptotic behavior of solutions $w(t,x)=\cT(t,0)w_0$, for $w_0\in\cA_f,\ t\rightarrow \pm\infty$.
See also \cite{furo24,roc25}.  
For example, eigenspaces of Floquet exponents for rotating waves $\cRP$, and of eigenvalues for homogeneous equilibria $\cEP$ and frozen waves $\cFP$\,, can be grouped by the zero numbers $z$ of their associated eigenfunctions. We obtain dimensions 2, for even $z>0$, and 1, for $z=0$.

More globally, and most importantly, zero number dropping \eqref{109} implies \emph{automatic transversality} of stable and unstable manifolds
\begin{equation}
\label{transv}
W^u(v_1)\transv W^s(v_2)
\end{equation}
for any pair of hyperbolic $v_1\,,\,v_2\in\cEP\cup\cFP\cup\cRP$.
See \cite{czro08, firowo04, hen85, roc25}.
In particular, $\bS^1$-equivariance then excludes homoclinic orbits $v_1=v_2$.
By decomposition \eqref{cAP}, heteroclinicity $\cHP$, i.e. $W^u(v_1)\cap W^s(v_2)\neq\emptyset$, is the only remaining option.
If all equilibria and rotating waves in the global attractor are hyperbolic, then transversality \eqref{transv} along all heteroclinic orbits $u{:}\ v_1\leadsto v_2$ and the celebrated $\lambda$-Lemma imply transitivity of heteroclinicity:
\begin{equation}
\label{transi}
v_1\leadsto v_2\leadsto v_3 \quad\Rightarrow\quad v_1\leadsto v_3\,.
\end{equation}
In absence of frozen waves, $\cFP=\emptyset$, the flow on the global attractor $\cAP$ becomes a Morse-Smale system.
In particular, the flow on $\cAP$ becomes \emph{structurally stable}: $C^2$-small perturbations of $f$ lead to $C^0$ orbit-equivalent flows on the global attractors.
See \cite{Palis, PalisSmale, PalisdeMelo} for general background on transversality and Morse-Smale systems.

In the hyperbolic case, compactness of $\cAP$ implies that $\cEP\cup\cFP\cup\cRP$ consists of finitely many $\bS^1$ orbits.
Passing to co-rotating coordinates \eqref{101c}, for some small $|c|\neq 0$, then eliminates all frozen waves, and Morse-Smale structural stability applies.
However, the perturbation of any frozen wave $v\in\cFP$ from $c=0$ to nonzero $c$ is not an orbit equivalence: a circle of equilibria is never equivalent to that same circle, as a nonstationary periodic orbit.
As a remedy, we could consider the induced semiflow, attractor, etc, on the orbit space $X/\bS^1$, which collapses $\bS^1$ orbits to single points.
Technical issues like transversality, the $\lambda$-Lemma, and structural stability would then have to be settled.

\sss{Connection graphs}\label{Conn}

\begin{figure}[t] 
\begin{center}
\includegraphics[width =  \textwidth]{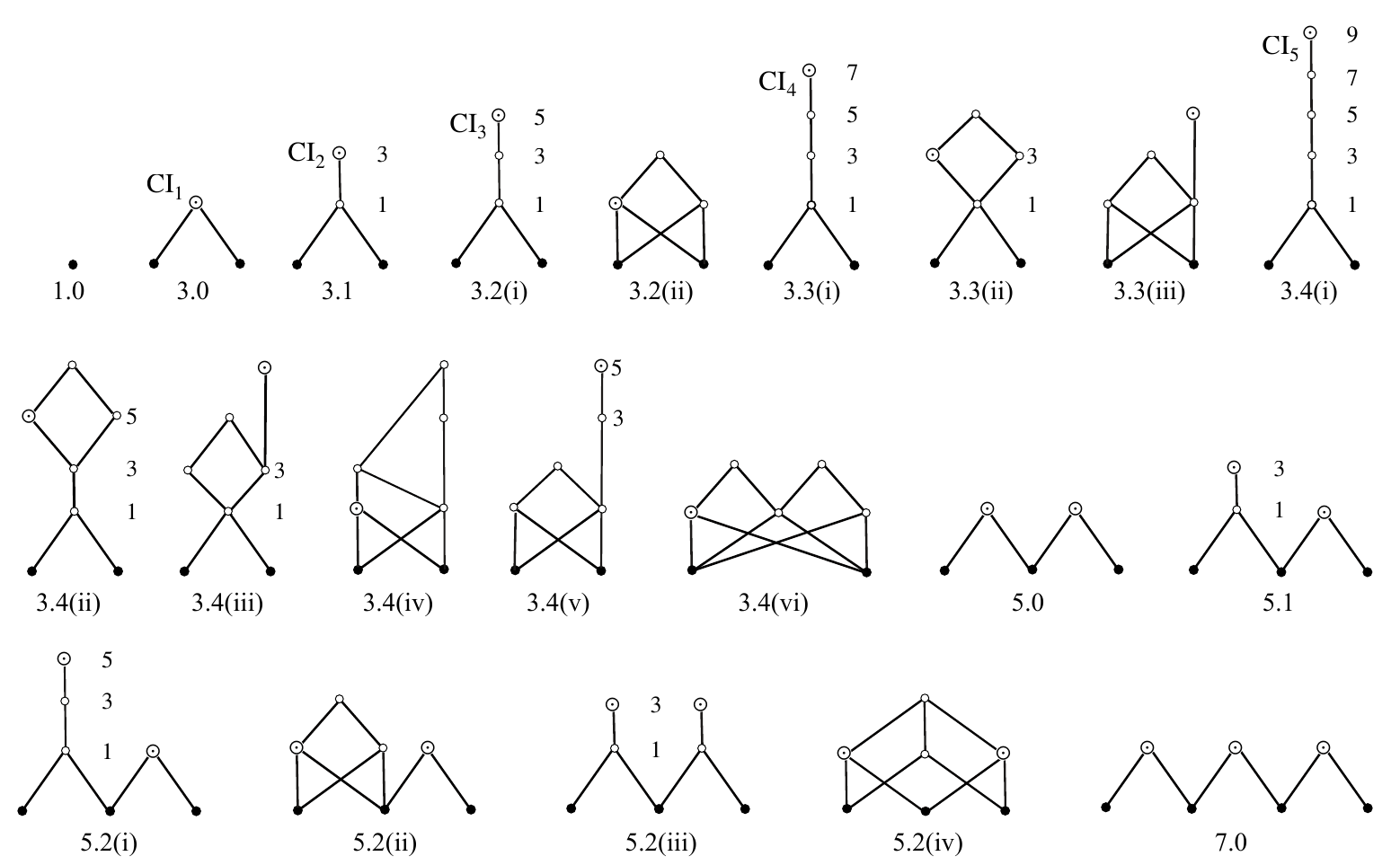}
\end{center}
\caption[Connection graphs for $n+q\leq7$]{\small\emph{ List of all 21 connection graphs $\cCP$ of dissipative PDEs \eqref{101} with $n$ homogeneous equilibria $e\in\cEP$ and $q$ frozen or rotating waves $v\in\cFP\cup\cRP$ (circles $\circ$), for $n+q\leq 7$.
The list is ordered by $n.q$, lexicographically.
The $(n+1)/2$ PDE stable homogeneous equilibria are marked by bullets $\bullet$.
Circled dots $\odot$ mark the $(n-1)/2$ PDE unstable homogeneous equilibria.
Vertical grading of each connection graph is by ascending Morse indices $i(e), i(v)$; see \eqref{ie}, \eqref{iv}.
Edges are oriented downhill, towards lower Morse index.
Nonadjacent jumps in Morse grading, only, are annotated explicitly.
For associated Sturm meanders see figure \ref{figmeanders}.
For our underlying classification by full lap signatures see table \ref{tab21}.
For the Chafee-Infante cases $\mathrm{CI}_{q+1}$ in 3.0, 3.1, 3.2(i), 3.3(i), 3.4(i) see section \ref{secCI}.
}}
\label{figgraphs}
\end{figure}

We use a more combinatorial coarsening, instead, based on directed \emph{connection graphs}.
The \emph{vertices} of the connection graph $\cCP$ are the $\bS^1$ orbits in $\cEP\cup\cFP\cup\cRP$.
By hyperbolicity assumption \eqref{HYP} and compactness of the global attractor $\cAP$, the number of vertices is finite.
Let $v_1\,,\,v_2$ denote any two vertices.
\emph{Directed edges} $v_1\leadsto v_2$ then signify the existence of any heteroclinic orbit $u$ from $v_1$ to $v_2$ in $\cHP$. 
We do not draw multiple edges to account for more than one such time orbit.
We call $f,g\in\SSOP$ and their global attractors $\cA_f^\cP,\cA_g^\cP$ \emph{connection equivalent}, when their directed connection graphs are isomorphic: $\cCP\cong\cC_g^\cP$.
Note how connection equivalence is a weaker notion, in general, than $C^0$ orbit equivalence of global attractors.
For recent surveys of various other notions, contexts, terminologies, and topologies related to connection graphs see \cite{Yorke-b, Yorke-a} and the many references there.
 
Since vertices are spatial $\bS^1$ orbits, we do not resolve in-phase convergence to rotating waves or specific fast convergence to frozen waves.
By \eqref{transi}, any connection graph $\cCP$ is transitive.
This allows us to replace $\cCP$ by the \emph{minimal graph} which generates the same transitive closure.
By transversality \eqref{transv} and $\bS^1$-equivariance, the directed connection graph is graded by the Morse indices $i$ from \eqref{ie}, \eqref{iv}:
\begin{equation}
\label{Ci}
v_1\leadsto v_2 \quad\Rightarrow\quad i(v_1)>i(v_2)\,,
\end{equation}
induces downhill orientation of all directed edges.
For illustration of the methods outlined in this survey, we enumerate all connection graphs $\cCP$ with up to seven vertices.

\begin{theorem}\label{th21}
Consider the parabolic PDE \eqref{101}, equivariant under $\bS^1\cong\bSO$, with Sturm nonlinearities $f\in\SSOP$ as defined in \eqref{sturmP}.
Let $n$ count the homogeneous equilibria $e\in\cEP$ and $q$ the frozen or rotating waves $v\in\cFP\cup\cRP$ in the global attractor $\cAP$.
Then the 21 connection graphs $\cCP$ of PDE \eqref{101} with $n+q\leq 7$ are enumerated in figure \ref{figgraphs}.
\end{theorem}

In \cite{firowo12b} we have already listed all 18 cases with $N:=n+2q\leq9$.
The cases were distilled by reduction to Neumann boundary values, and a subsequent brute force scan of the 43 Sturm permutations known from \cite{fiedlertatra,firo96}.
Three of the remaining 18 cases exceeded the current bound $n+q\leq7$.
Therefore, only the six cases 3.4(i)--3.4(vi) are new in the present paper.
In section \ref{Nlap}, however, we present an enumerative method, based on period maps \eqref{TlP} and lap signatures, which aims for the setting of theorem \ref{th21}, directly.
See table \ref{tab21} below.
Moreover, the new approach offers rare conceptual insight into the little understood link between rotating wave ODEs like \eqref{103} and the global dynamics of PDE \eqref{101}.

\sss{Outline}\label{Out}

A more detailed outline of our survey is as follows.
Theorem \ref{th21} will be proved in section \ref{PfConn}, based on a general method which we survey in sections \ref{PNbc}--\ref{Nlap}.
In section \ref{PNbc}, we start with a homotopy from the $\bSO$ Sturm class $f\in\SSOP$ to the special case $f=g$ of spatially reversible nonlinearities $g\in\SO$.
The resulting enhanced $\bO$-equivariance with respect to spatial shifts and flips $x\mapsto 2\pi-x$ eliminates all rotating waves, $\cR_g^\cP=\emptyset$, and makes all frozen waves $v\in\cF_g^\cP$ reflection-symmetric.
In particular, suitably shifted copies of $v$ become equilibria $v\in\cFN$ of the associated Neumann problem on the half-interval $0<x<\pi$.
Moreover, we show how closely the connection graphs $\cC_g^\cP$ and $\cCN$ are related.
Section \ref{Nmeander} recalls the combinatorial characterization of Neumann Sturm global attractors $\cAN$ by certain meander permutations $\sigma_g$. 
Meander permutations arise, for example, from a shooting approach to equilibria, i.e. to solutions $v(x)$ of the Neumann boundary value problem for ODE \eqref{102}.
By integrability of \eqref{102}, in the spatially reversible case $f=g\in\SO$, we can alternatively classify Sturm global attractors $\cAN$ by the full lap signatures $\mathfrak{S}(T_g)$ of the period map $T_g$ associated to the nonlinearity $g$; see section \ref{Nlap}.
In section \ref{PfConn} we prove theorem \ref{th21}, by enumerating all 21 full lap signatures of period maps, for $n+q\le 7$. 
In addition, we list all associated stylized shooting meanders. 
Section \ref{Delay} reviews some related results for delay differential equations (DDEs) \eqref{501} with positive monotone feedback.
In particular, we revisit the three-dimensional \emph{Vas tulip} \cite{vas17} as an example of a global attractor $\cA_T$ with $n=5$ equilibria and $q=4$ rotating waves; see figure \ref{figtulip}.
Surprisingly, a frozen version of the Vas tulip also appears among the global Sturm attractors in the spatially reversible class of $g\in\SO$.
See the construction via $g\in\SO$ and its meander, in figure  \ref{figTNM}.
The associated parabolic connection graphs are derived in figure \ref{figNP}.
Although the semiflows of monotone feedback DDEs share a discrete Lyapunov functional like the zero number \eqref{z}, \eqref{109} with PDEs \eqref{101}, the precise details of this fascinating correspondence remain to be explored further.
See in particular \cite{furo24}.

\section{Periodic versus Neumann boundary conditions}\label{PNbc}

In this section we compare the connection graphs $\cC_{f\,,\,g}^{\cP,\,\cN}$ of the global attractors $\cA_{f\,,\,g}^{\cP,\,\cN}$ under periodic boundary conditions  \eqref{101} versus the Neumann case \eqref{PDEg}. 
More precisely, we reduce the periodic $\bSO$-case $f\in\SSOP$ to the spatially reversible Neumann case \eqref{grev}, i.e.~to $g\in\SO$ with hidden $\bO$-symmetry.
See definition \ref{defsturm} for notation.

The comparison will proceed in two steps.
First, we compare periodic and Neumann boundary conditions for the special case that $f=g\in\SO$ already satisfy spatial reversibility \eqref{grev}.
In other words, we will compare the global attractors $\cgAP,\, \cAN$ and their connection graphs $\cgCP,\, \cCN$.
In a second step, we then recall how to reduce the general $\bSO$-case $f\in\SSOP$ to the spatially reversible $\bO$-case $f=g\in\SO$, without affecting the connection graph $\cCP$.
Our exposition follows the account of \cite{firowo04} in \cite{firowo12b}.

\sss{Spatial reversibility}\label{Prev2Nrev}

For our first step, $f=g$, we recall that the Neumann interval $0<x<\pi$ only spans half the circle $0<x<2\pi$ of $x\in\bS^1=\R/2\pi\Z$.
Under spatial reversibility \eqref{grev}, the Neumann case corresponds to the dynamics of the periodic case in its invariant (isotropy) subspace of flip-invariant solutions $\kappa u=u$; see \eqref{flip}.
Hyperbolicity of equilibria is also preserved, but the precise Morse indices may change.

The separated Neumann boundary conditions make \eqref{PDEg} gradient-like \cite{firaro, lafi, mat88}.
This eliminates time-periodic orbits, and the global attractor $\cAN$ of \eqref{PDEg} decomposes as
\begin{equation}
\label{cAN}
\cAN=\cEN\cup\cFN\cup\cHN.
\end{equation}
As before, $\cEN$ denotes the spatially homogeneous equilibria which are not affected by periodic versus Neumann boundary conditions, and $\cFN$ lists the remaining spatially nonhomogeneous equilibria:
\begin{equation}
\label{EFN}
\cEN=\{e_1,\ldots,e_n\} \quad\textrm{and}\quad \cFN=\{v_1,\ldots,v_{2q}\}\,,
\end{equation}
totaling $N=n+2q$ equilibria.
Indeed, $N$ and $n$ are both odd, by dissipativity.
Again $\cHN$ collects the heteroclinic orbits between any equilibria.

Let us explore the consequences of spatial reversibility \eqref{grev} on equilibria, next.
Any PDE equilibrium $v\in \cEN\cap\cFN$ solves the second order ODE \eqref{102}, for $f:=g$, i.e.
\begin{equation}
\label{vpg}
\begin{aligned}
v_x&=p\,,\\
p_x&=-g(v,p)\,.
\end{aligned}
\end{equation}
The nonhomogeneous PDE equilibria $v\in\cFN$ become nonstationary periodic ODE orbits $(v,p)$.
Even in absence of reversibility \eqref{grev}, we have seen how any nonstationary periodic ODE orbit $(v(x),p(x))$ intersects the horizontal axis twice: at its minimum $\underline{v}$ and at its maximum $\overline{v}$; see \eqref{minmax}.
Reversibility \eqref{grev}, however, makes the periodic ODE orbit reflection symmetric.
In particular it spends half its minimal period $T_g(v)>0$ above the horizontal axis, and the other half below; compare \eqref{Tpm} with $p_-(v)=-p_+(v)$ due to spatial reversibility.
Hence the extrema $v_x=0$ of the periodic ODE orbit $x\mapsto v(x)$ are spaced half a period $T_g(v_0)$ apart.
By suitable phase shifts of $v(x)$, each nonstationary periodic ODE orbit $(v(x),p(x))$ of \eqref{vpg} with minimal period
\begin{equation}
\label{TlN}
0<T_g(v)=2\pi/\ell\,, \quad \ell\in\N
\end{equation}
therefore gives rise to a \emph{pair of Neumann equilibria}, bijectively:
\begin{equation}
\label{vpair}
\cF_g^\cP\ni v \ \longleftrightarrow \  \underline{v},\, \overline{v} \in\cFN \,.  
\end{equation}
We therefore call $v\in\cFN$ \emph{frozen equilibria}.
Indeed, we have the distinct alternatives $\underline{v}(0)=\min v$ and $\overline{v}(0)=\max v$.
We also note that the integer fraction
\begin{equation}
\label{lap}
\ell=\ell(v)=z^\cN(v_x)+1\in \N
\end{equation}
in \eqref{TlN} coincides with the Matano \emph{lap number}; see \cite{mat82}.
Note $\ell(\underline{v})=\ell(\overline{v})$.
The lap and zero numbers $\ell,\,z^\cN$ refer to the Neumann half interval $0<x<\pi$;
but compare also \eqref{TlP}.
On bounded intervals, the zero number $z^\cN$ of \eqref{z} only counts strict interior sign changes; therefore the Neumann boundary values $v_x=0$ at $x=0,\pi$ do not contribute to the lap number.

We may also call the min/max pair correspondence \eqref{vpair} a \emph{hidden symmetry}, in the sense of \cite{gomasc84}.
Indeed $\bO$-equivariance of \eqref{101} under periodic boundary conditions, for $f=g\in\SO$, represents the $\bSO$ orbit of the frozen wave equilibrium $v$ by the two $x$-shifted copies $\underline{v},\, \overline{v}$, each invariant under the flip $x\mapsto -x$ of \eqref{flip}.
In fact, each $x$-periodic extension $\underline{v},\, \overline{v}$ is fixed under the dihedral isotropy subgroup
\begin{equation}
\label{isotropy}
\mathbb{D}_\ell:=\Z_\ell\rtimes\langle\kappa\rangle \leq \bO\,.
\end{equation}
Here again, $\ell$ denotes the (Neumann) Matano lap number \eqref{lap}.
For later reference note
\begin{equation}
\label{lparity}
\begin{aligned}
   \underline{v}(\pi)&=\overline{v}(0) \quad \textrm{and}\quad \overline{v}(\pi)=\underline{v}(0), \textrm{ for odd } \ell\,, \textrm{ and}\\
\underline{v}(\pi)&=\underline{v}(0) \quad \textrm{and}\quad  \overline{v}(\pi)=\overline{v}(0), \textrm{ for even } \ell\,.
\end{aligned}
\end{equation}
In this language, we can now summarize the results already established in \cite{firowo12b}.

\begin{theorem}\label{thgCPN}
Assume $g\in\SO$ is spatially reversible; see \eqref{grev} and definition \ref{defsturm}.
Then the quite distinct global attractors $\cA_g^\cP$ and $\cAN$, under periodic and Neumann boundary conditions respectively, share the same connection graphs, up to hidden symmetry.
More precisely, 
\begin{equation}
\label{gCPN}
\cC_g^\cP\ \cong\ \cCN /\sim\,.
\end{equation}
Here $\sim$ identifies vertex pairs $\underline{v},\overline{v}\in\cFN$ of frozen Neumann equilibria which share the same extremal values, i.e. which are related by the hidden $\bO$-symmetry \eqref{vpair}.
\end{theorem}

\sss{Towards spatial reversibility}\label{P2Prev}

In our second step, it remains to perform a reduction from the general $\bSO$-case $f\in\SSOP$ to the spatially reversible $\bO$-case $f=g\in\SO$, again without affecting the connection graph $\cCP$. See \cite{firowo04} for details, which prove the following theorem.

\begin{theorem}\label{thCPtau}
Assume $f\in\SSOP$ is not necessarily $x$-reversible.
Then there exists a homotopy $f^\tau\in\SSOP,\ 0\leq\tau\leq1$, of uniformly dissipative nonlinearities in the same class, such that $f=f^0$, but $g:=f^1\in\SO$ satisfies the additional spatial reversibility condition \eqref{grev}.
The homotopy fixes homogeneous equilibria.
Hyperbolicity  \eqref{HYP} tracks rotating and frozen waves.
In particular, the connection graphs remain isomorphic throughout the homotopy:
\begin{equation}
\label{CPtau}
\cCP\cong\cC_{f^\tau}^\cP\cong\cgCP\,.
\end{equation}
\end{theorem}

Combining both theorems, \eqref{CPtau} and \eqref{gCPN} provide a Neumann description of the $\bSO$ connection graph $\cCP$ for all, not necessarily reversible, nonlinearities $f$.

\begin{corollary}\label{corCPN}
For any nonlinearity $f\in\SSOP$ there exists a spatially reversible nonlinearity $g\in\SO$ such that
\begin{equation}
\label{CPN}
\cCP\cong\cCN /\sim\,.
\end{equation}
As in \eqref{vpair} and theorem \ref{thgCPN}, $\underline{v}\sim\overline{v}$ identifies vertices by the hidden $\bO$-symmetry \eqref{vpair}.
\end{corollary}
We conclude this section with a sketch of the homotopy $f^\tau\in\SSOP$ used in theorem \ref{thCPtau}.
The homotopy consists of two different parts, for $0\leq\tau\leq1/2$ and for $1/2\leq\tau\leq 1$.
See sections \ref{Freeze} and \ref{Symm}.

\sss{Freeze \ldots}\label{Freeze}

The first homotopy, for $0\leq\tau\leq1/2$, freezes all rotating waves; see section 4 of \cite{firowo04} for complete details.
We recall that rotating waves $\mathbf{v}(x):=(v(x),p(x))$ are the solutions of the second order ODE \eqref{103} with (not necessarily minimal) period $2\pi$.
In other words,
\begin{equation}
\label{vp}
\begin{aligned}
v_x&=p\,,\\
p_x&=-f-cp\,,
\end{aligned}
\end{equation}
with $f=f(v,p)$.
In \cite{mana97}, Matano and Nakamura discovered that the nonstationary periodic orbits $\mathbf{v}$ of \eqref{vp}, for whatever fixed speed parameters $c\in\R$, foliate an invariant open region $\mathbf{C}\subset\R^2$ which we call the \emph{cyclicity set}.
In other words,
\begin{equation}
\label{cycset}
\mathbf{C} := \bigl\{\,\mathbf{v}_0=(v_0,p_0)\in\R^2\ \big|\ x\mapsto \mathbf{v}(x) \textrm{ is nonstationary periodic, for some } c\in\R\, \bigr\}\,.
\end{equation}
Remarkably, the associated wave speed $c$ is unique, for each periodic point$(v_0,p_0)\in\mathbf{C}$; see \cite{mana97}.
Indeed, suppose indirectly that $\mathbf{v}_1\,,\mathbf{v}_2$ are two rotating or frozen waves at different wave speeds $c_1>c_2$\,.
Assume any intersection, at $\mathbf{v}_1(x_0)=(v_0,p_0)=\mathbf{v}_2(x_0)$.
Abbreviating the right side of ODEs \eqref{vp} at the intersections by $\mathbf{f}_j$\,, we obtain the determinant, cross product, or alternating form
\begin{equation}
\label{^}
\mathbf{f}_1\wedge\mathbf{f}_2\ =\ (c_1-c_2)\,p_0^2\ >\ 0\,,
\end{equation}
for $p_0\neq0$.
This shows that $\mathbf{v}_2$ can only cross $\mathbf{v}_1$ from right to left, at $p_0\neq0$.
Two more differentiations with respect to $x$ extend this observation to $p_0=0$.
Such uni-directional crossing, however, contradicts periodicity of $\mathbf{v}_2$ and the Jordan curve theorem for $\mathbf{v}_1$\,.

\begin{figure}[t] 
\begin{center}
\includegraphics[width = 0.7\textwidth]{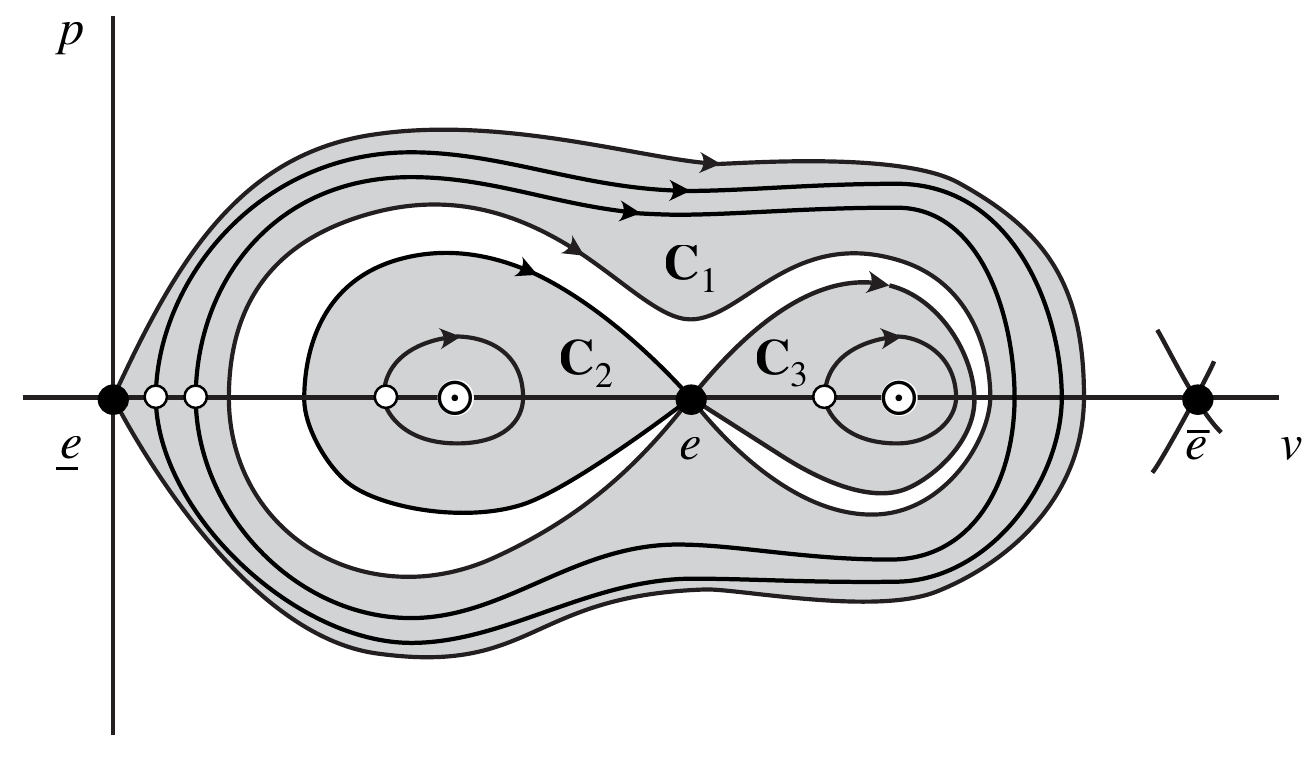}
\end{center}
\caption[Cyclicity set of the Vas tulip]{\small\emph{ Example of an open cyclicity set $\mathbf{C}=\mathbf{C}_1\cup\mathbf{C}_2\cup\mathbf{C}_3$ (shaded gray) of ODE \eqref{vp}.
The connected components $\mathbf{C}_2$ and $\mathbf{C}_3$ are punctured disks.
The connected component $\mathbf{C}_1$ is an annulus. 
Note the $n=5$ homogeneous equilibria, alternatingly of saddle and center type, $\bullet$ and $\odot$.
The centers are the two punctures.
The middle saddle is annotated by $e$.
The annulus $\mathbf{C}_1$ contains two orbits of period $2\pi$.
The punctured disks $\mathbf{C}_2, \mathbf{C}_3$ contain one such periodic orbit, each.
These periods are minimal.
Therefore we obtain $q=4$ rotating/frozen waves.
See also the full lap signature \eqref{304}, and the Vas tulips of figures \ref{figtulip} and \ref{figTNM}(left).
}}
\label{figcyc}
\end{figure}

The connected components $\mathbf{C}_k$ of the cyclicity set $\mathbf{C}$ consist of \emph{punctured disks} and \emph{annuli}. 
Each puncture is a homogeneous equilibrium $\mathbf{e}_0=(e_0,0)$ with $f_u(e_0,0)>0$, which becomes a linearized \emph{center} for limiting wave speed $c=-f_p(e_0,0)$.
By an implicit function theorem along the horizontal $v$-axis, the surrounding punctured disks extend outward until the upper or lower return times $T_\pm$ to the horizontal axis become unbounded; see \eqref{Tpm}.
In fact, the periodic orbits would have to approach some ODE \emph{saddle} equilibrium $\mathbf{e}_1=(e_1,0)$ with $f_u(e_1,0)<0$, for some finite limiting wave speed $c_*\in\R$\,.
Indeed, $c_*=\pm\infty$ requires convergence to the horizontal axis and contradicts periodicity; see also \cite{mana97} for an a priori bound.
For the same reason, the remaining annular cyclicity components experience period \emph{blow-up} $T_\pm\rightarrow\infty$, and possess ODE saddles at their outer and inner boundaries.
See figure \ref{figcyc} for an example of a cyclicity set $\mathbf{C}$ with a single annular component $\mathbf{C}_1$ which surrounds two punctured disks $\mathbf{C}_2\,,\,\mathbf{C}_3$\,.
Along \emph{blow-up boundaries}, we encounter homoclinic or heteroclinic loops involving one or several saddles.
The precise planar configurations of these boundaries of the cyclicity set is neither unique, not even for fixed integrable Sturm involutions $\sigma=\sigma_g$, nor relevant for the resulting connection graphs.
Indeed \eqref{TlP} reminds us to be concerned with $T_++T_-=T\leq 2\pi$, only.

In summary, the unique wave speeds $c$ and the minimal periods $T>0$ define $C^2$ functions
\begin{equation}
\label{cT}
c,T:\quad \mathbf{C}\rightarrow \R
\end{equation}
on the cyclicity set $\mathbf{C}$.
Note how the ODE-based maps $c,T$ are defined, independently of periodic or Neumann PDE boundary conditions $\cP$ or $\cN$.

For a moment, suppose $c$ were defined globally.
Then we can define a \emph{freezing homotopy}
\begin{equation}
\label{f1}
f^\tau(v,p):=f(v,p)+2\tau c(v,p)p
\end{equation}
in $\SSOP$.
Indeed, the homotopy preserves all homogeneous equilibria, frozen and rotating waves, including shapes and spatial periods.
Wave speeds $c^0$ are slowed to $c^\tau=(1-2\tau)c^0$.
At $\tau=1/2$, this freezes all rotating waves.
Remarkably, freezing also preserves hyperbolicity; compare \eqref{101c}.
By automatic transversality \eqref{transv}, it also preserves the connection graph along our first homotopy $0\leq\tau\leq1/2$:
\begin{equation}
\label{CPtau1}
\cCP\cong\cC_{f^\tau}^\cP\cong\cC_{f^{1/2}}^\cP\,.
\end{equation}
At $\tau=1/2$, the resulting ODE \eqref{vp} becomes integrable on the cyclicity region $\mathbf{C}$.
Indeed, the wave speeds $c=c(v,p)$ provide a first integral, with the periodic orbits \eqref{cycset} as level curves.
In other words, $c$ solves the \emph{nonlinear hyperbolic conservation law}
\begin{equation}
\label{chyp}
(pc+f)c_p-pc_v=0\,.
\end{equation}
The conservation law degenerates at the ODE equilibria $p=f=0$, i.e.~at saddles and centers $(e,0)$.
It is therefore not surprising when the cyclicity set $\mathbf{C}$ contains gaps, globally, where $c$ cannot be defined consistently.
But global consistency is not required for our freezing homotopy, either.
Indeed it is sufficient to define freezing on the subset of $\mathbf{C}$ with minimal periods $T\leq 2\pi$, to capture all rotating and frozen waves relevant for our domain of $x\in\bS^1=\R/2\pi\Z$.
It is then easy to $C^2$ extend $c$ globally, to all of $\R^2$, without introducing new spurious periodic orbits of minimal periods $T\leq 2\pi$.
This establishes the connection graph equivalence \eqref{CPtau1}.

\sss{\ldots and symmetrize}\label{Symm}

To achieve spatial reversibility \eqref{grev}, we sketch how a second homotopy $1/2\leq\tau\leq1$ symmetrizes $f^{1/2}$ and all its frozen waves; see section 5 of \cite{firowo04} for full details.
We proceed via the underlying ODEs \eqref{103}, alias \eqref{vp}, with $f=f^\tau$, but now with frozen $c=0$.
Again it is sufficient to perform the homotopy on the subset of the integrable cyclicity set with minimal periods $T\leq 2\pi$.
We perform the homotopy on nonstationary periodic orbits $\mathbf{v}=(v,p)$ which foliate the cyclicity set at $\tau=1/2$.
Fix one such orbit.
It only crosses the horizontal $v$-axis $p=0$ at $\underline{v}=\min v$ and $\overline{v}=\max v$; see \eqref{minmax}.
Any other value $v$ along the orbit is crossed twice, at unique values $p_-(v)<0<p_+(v)$ of $p$.
The homotopy $f^\tau$ for $1/2\leq\tau\leq1$ is then defined by \emph{harmonic homotopy} of our periodic orbit, i.e.
\begin{equation}
\label{ptau}
1/p_\pm^\tau := 2 (1-\tau)\cdot 1/p_\pm + (\tau-\tfrac{1}{2})\cdot (1/p_\pm-1/p_\mp)\neq 0\,,
\end{equation}
for any non-extremal value $v$.
Note $p_\pm^{1/2}=p_\pm$ and the symmetrization $1/p_+^1=-1/p_-^1$\,.
For second order ODEs \eqref{102}, alias \eqref{vp} with $c=0$, orbits $\mathbf{v}=(v,p(v))$ define nonlinearities $f=f(u,v)$ along them.
Indeed,
\begin{equation}
\label{fvp}
f(v,p)=-p\,dp/dv\,.
\end{equation}
The homotopy $p=p_\pm^\tau(v)$ of periodic orbits therefore defines a homotopy $f^\tau$ of nonlinearities, in our subset $T\leq 2\pi$ of the cyclicity region.
Symmetry of periodic orbits, at $\tau=1$, implies spatial reversibility \eqref{grev} of $g:=f^1$.
Moreover, minimal periods $T$ and hyperbolicity remain preserved by harmonic homotopy.
Indeed, \eqref{Tpm} and \eqref{ptau} imply
\begin{equation}
\label{eqT}
T=\int_{\underline{v}}^{\overline{v}} \,(1/p_+^\tau-1/p_-^\tau)\,dv\,.
\end{equation}
Again we refer to \cite{firowo12b} for further details, including hyperbolicity and $C^2$ regularity.
As before, automatic transversality \eqref{transv} then preserves the connection graph along our second homotopy $1/2\leq\tau\leq1$ of $f^\tau\in\SSOP$:
\begin{equation}
\label{CPtau2}
\cC_{f^{1/2}}^\cP\cong\cC_{f^\tau}^\cP\cong\cgCP\,.
\end{equation}
Combining the two homotopies \eqref{CPtau1} and \eqref{CPtau2} establishes claim \eqref{CPtau} of theorem \ref{thCPtau}.

\section{Neumann attractors and Sturm permutations}\label{Nmeander}

Corollary \ref{corCPN} above has reduced the periodic connection graph $\cCP$ of any $f\in\SSOP$ to the Neumann connection graph $\cCN$ for some spatially reversible $g\in\SO$ on the half-interval $0\leq x\leq \pi$.
The hidden $\bO$-symmetry of the Neumann problem has led to the equivalence relation $\sim$ in \eqref{gCPN}, which merges min/max vertex pairs $\underline{v},\overline{v}\in\cFN$ of spatially nonhomogeneous Neumann equilibria in \eqref{vpair} sharing the same extrema.

Section \ref{Meander} recalls our results on connection graphs $\cCN$ of Neumann problems.
We admit general, possibly $x$-dependent $g\in\SxN$; see \eqref{sturmNx}.
This part is based on \cite{firo96, firo99, firo00, furo91}.
See also \cite{firo14, firo13, firo22,rofi24} for the topology of the \emph{Thom-Smale complex} of heteroclinic orbits, alias unstable manifolds, as a regular topological cell complex,
and \cite{kar17} for beautiful illustrations and further mathematical connections.
Following \cite{firowo12a}, section \ref{IntInv} then surveys the modifications of $\cCN$ which arise in the restricted reversible case $g\in\SO$.

\sss{Meander permutations and Sturm permutations}\label{Meander}

We recall some basic definitions.
A \emph{meander} $\cM$ is an oriented planar Jordan curve, infinitely extended and running asymptotically from southwest to northeast, which transversely intersects an oriented horizontal 
line at a finite number $N$ of points \cite{arn89}.
Labeling intersections along $\cM$, their sequence ordered along  the horizontal axis defines a permutation $\sigma=\sigma_\cM\in S(N)$ which we call a \emph{meander permutation}.

But let us return to the Sturm global attractor  $\cAN$ of the general Neumann problem \eqref{PDEg}. 
A characterization of this global attractor in combinatorial terms is provided by a permutation defined on the set of (homogeneous or nonhomogeneous) Neumann equilibria $v_1,\ldots,v_N\in\cE_g\cup\cF_g^\cN$. 
We label equilibria by their ordering at $x=0$:
\begin{equation} \label{201}
\underline{e}=v^\cN_1(0) < v^\cN_2(0) < \dots < v^\cN_N(0)=\overline{e} \,.
\end{equation}
The ordering of the same equilibria at the other Neumann boundary $x=\pi$ will be different, in general.
This defines a \emph{Sturm permutation} $\sigma=\sigma_g^\cN$ such that
\begin{equation} \label{202}
\underline{e}=v^\cN_{\sigma(1)}(\pi) < v^\cN_{\sigma(2)}(\pi) < \dots < v^\cN_{\sigma(N)}(\pi)=\overline{e}\,.
\end{equation}
In particular $\sigma_g^\cN(1)=1$ and $\sigma_g^\cN(N)=N$; see also theorem \ref{th1}(ii) below.

Let $\mathbf{v}(x;\mathbf{v}_0)$ denote solutions $\mathbf{v}=(v,p)$ of second order equations \eqref{102}, \eqref{vp}, frozen at $c=0$, with initial condition $\mathbf{v}=\mathbf{v}_0$ at $x=0$.
A shooting approach to the Neumann boundary value problem then studies intersections of the shooting curve
\begin{equation}
\label{shoot}
\cM:\quad v_0\mapsto \mathbf{v}(\pi;v_0,0)
\end{equation}
with the horizontal Neumann axis $p=0$.
By standard Sturm-Liouville theory, transverse intersections are equivalent to PDE hyperbolicity of the equilibrium $v(x)$.
Thus the hyperbolic case corresponds to transverse meanders $\cM$, and the Sturm permutation $\sigma_g=\sigma_\cM$ defined by \eqref{201}--\eqref{202} is a meander permutation. 
Conversely, we often represent any meander permutation $\sigma$ by a \emph{stylized meander}, which consists of alternating upper and lower arcs over horizontal positions $\sigma^{-1}(j),\ \sigma^{-1}(j+1)$.
See figure \ref{figmeanders} for 21 examples, and \cite{kar17} for many more illustrations.

The main results of \cite{firo96,firo00} can now be summarized as follows.

\begin{theorem}\label{thfCAN}
Consider general nonlinearities $g\in\SxN$ for the Neumann problem \eqref{PDEg}; see definition \ref{defsturm}.
Then the associated Sturm permutation $\sigma_g^\cN\in S(N)$ is a meander permutation and determines
\begin{enumerate}[(i)]
  \item all Morse indices $i(v^\cN_j)$ and all zero numbers $z(v^\cN_j-v^\cN_{j'})$, for $1\leq j\neq j'\leq N$;
  \item the connection graph $\cC_g^\cN$; and
  \item the global attractor $\cA_g^\cN$, up to $C^0$ orbit equivalence.
\end{enumerate}
\end{theorem}

The formul\ae\ for Morse indices and zero numbers are explicit, and allow an explicit construction of the connection graph. 
Abbreviating $\sigma_{kj}:= \sign\,(\sigma^{-1}(k)-\sigma^{-1}(j))$,
\begin{equation}
\label{ij}
\begin{aligned}
   i^\cN(v_j^\cN)=i_j&:=\sum_{k=1}^{j-1}\,(-1)^{k+1}\sigma_{k+1,k}\,, \\
   z^\cN(v^\cN_j-v^\cN_{j'})=z_{jj'} &  := i_j+\tfrac{1}{2}\big((-1)^{j'}\sigma_{j'j}-1\big)+\sum_{j<k<j'}(-1)^k\sigma_{kj}\,,
\end{aligned}
\end{equation}
for $1\leq j<j'\leq N$.
We also fix $i_1=i_N=0$ and symmetrically extend $z_{j'j}:=z_{jj'}$\,.

Based on $i_j$ and $z_{jj'}$ the edges of the connection graph $\mathcal{C}^{\mathcal{N}}_{g}$, viz.~the heteroclinic orbits, can be determined as follows; see \cite{wol02} and the appendix in \cite{firo19}.
Zero number dropping \eqref{109} \emph{blocks}, i.e. prevents, Neumann heteroclinic orbits $v_j\leadsto v_{j'}$ if there exists another homogeneous or frozen Neumann equilibrium $v_k$\,, between $v_j$ and $v_{j'}$ at $x=0$, such that 
\begin{equation}
\label{block}
z_{jk}=z_{jj'}=z_{j'k}\,.
\end{equation}
We call the equilibria $v_j\,,\,v_{j'}\ \cN$-\emph{adjacent} if there does not exist any blocking equilibrium.
The explicit version of theorem \ref{thfCAN}(ii) then reads: $v_j\leadsto v_{j'}$, if and only if  $i_j>i_{j'}$ and $v_j\,,\,v_{j'}$ are $\cN$-adjacent.
The proof of claim (iii) still requires a homotopy between $f$ and $g$ in a spatially discretized setting \cite{firo00}.

The powerful reduction of global PDE dynamics to the discrete math of certain meander permutations, above, renders it centrally important to determine the range of all Sturm permutations $\sigma_g^\cN$ in $S(N)$, abstractly.
See \cite{firo99} for details.
Shooting \eqref{shoot} identifies $\sigma_g^\cN$ as meandric.
The formal Morse numbers $i_j$\,, defined by \eqref{ij} for any $\sigma\in S(N)$, have to become Morse indices  $i_j=i(v_j^\cN)$ in case $\sigma=\sigma_g^\cN$, and are therefore constrained to be nonnegative.
In general, we therefore call a permutation $\sigma\in S(N)$ \emph{Morse}, if $i_j\geq0$ holds, for all $j$.

\begin{theorem}\label{th1}
A meander permutation $\sigma\in S(N)$ is  a \emph{Sturm permutation}, 
i.e. $\sigma=\sigma^\cN_g$ for some $g\in\SxN$, if and only if
\begin{enumerate}[(i)]
  \item $N$ is odd;
  \item $\sigma$ is \emph{dissipative}, i.e.~$\sigma(1)=1$ and $\sigma(N)=N$; and
  \item $\sigma$ is Morse.
\end{enumerate}
\end{theorem}

\sss{Integrable Sturm involutions}\label{IntInv}

We now attend to the $\bO$-case of nonlinearities $g\in\SO$ which satisfy the additional spatial reversibility condition \eqref{grev}.
We mainly follow \cite{firowo12a}.
See also \cite{firowo04, firowo12b, roc07, roc25, rofi14} for further details.

Expectedly, spatial reversibility further restricts the range of Sturm permutations $\sigma=\sigma_g^\cN$ which can be realized in this class.
To explore these restrictions, we first collect some further properties of such $\sigma=\sigma_g^\cN\in S(N)$.

The frozen nonhomogeneous min/max pairs $v=\underline{v}_j,\overline{v}_j\in\cFN$ of \eqref{vpair} come in two flavors, depending on the parity of their lap number $\ell(v)$ in \eqref{lap}. 
By \eqref{lparity}, each frozen min/max pair $\underline{v}_{\,\alpha},\overline{v}_\alpha\in\cFN$ with \emph{odd lap number} $\ell$ contributes a \emph{2-cycle} 
\begin{equation}
\label{2cycle}
(\underline c_{\,\alpha} \ \overline c_\alpha)\,
\end{equation}
to the Sturm permutation $\sigma=\sigma_g^\cN$.
Without loss, we sort the 2-cycle as $\underline c_{\,\alpha} < \overline c_\alpha$\,.
The remaining boundary values, of homogeneous equilibria and of frozen min/max pairs with \emph{even lap number} $\ell$, contribute fixed points.
In particular, $\sigma$ is an \emph{involution},
\begin{equation} \label{203}
\sigma\circ\sigma = \id\,.
\end{equation}
Specifically, $\sigma$ is the finite product of disjoint 2-cycles \eqref{2cycle} which arise from odd lap numbers $\ell(\underline{v}_{\,\alpha})=\ell(\overline{v}_\alpha)$.

The following terminology is natural, when we recall how \eqref{vpair} associates integer 2-cycles \eqref{2cycle} with periodic ODE orbits \eqref{vpg} of frozen equilibria:
\begin{itemize}
\item two distinct 2-cycles $(\underline c_{\,\alpha} \ \overline c_\alpha)$ and 
$(\underline c_{\,\beta} \ \overline c_\beta)$ are called \emph{intersecting}, if the corresponding intervals have a 
nonempty intersection, i.e. 
$(\underline c_{\,\alpha}, \overline c_\alpha) \cap (\underline c_{\,\beta}, \overline c_\beta) \ne \emptyset$\,;
\item intersecting 2-cycles $(\underline c_{\,\alpha} \ \overline c_\alpha)$ and 
$(\underline c_{\,\beta} \ \overline c_\beta)$ are called \emph{nested} if one of the intervals strictly contains 
the other, i.e. $(\underline c_{\,\beta} - \underline c_{\,\alpha})(\overline c_\alpha - \overline c_\beta)>0$\,;
\item nested 2-cycles $(\underline c_{\,\alpha} \ \overline c_\alpha)$ and 
$(\underline c_{\,\beta} \ \overline c_\beta)$ are \emph{centered} if the mid-points of the respective 
intervals coincide, i.e. 
$\underline c_{\,\beta} - \underline c_{\,\alpha} = \overline c_\alpha - \overline c_\beta$\,.
\end{itemize}

A fixed point $e$ of $\sigma$ is called $\sigma$\emph{-stable}, if it possesses formal Morse index $i_e=0$; see \eqref{ij}.
For Sturmian $\sigma=\sigma_g$, we recall how $i_e=0$ implies that $e$ is an ODE saddle; see \eqref{ie}.
The notion of $\sigma$-stability is completely abstract, however, if we recall the formal Morse numbers $i_j$ associated to any permutation $\sigma\in S(N)$.
Similarly, we define the $\sigma$\emph{-stable core} $C_\alpha$ of a 2-cycle $(\underline c_{\,\alpha} \ \overline c_\alpha)$ to consist of its interior $\sigma$-stable points. 
We say that two nested 2-cycles of $\sigma$ are \emph{core-equivalent} 
if both share the same $\sigma$-stable core. 
The above terminology cumulates in the following definition.

\begin{definition}\label{defint}
An involution $\sigma$ is called \emph{integrable} if the following three properties hold:
\begin{enumerate}[(i)]
\item any two intersecting 2-cycles are nested;
\item any two core-equivalent 2-cycles are centered; and
\item any pair of non-nested 2-cycles is separated by a $\sigma$-stable point.
\end{enumerate}
\end{definition}

The necessary conditions collected in this definition characterize the additional restrictions of Sturm permutations $\sigma_g^\cN$ in the spatially reversible class; see \cite{firowo12a, roc25, rofi14} for full details.

\begin{theorem}\label{th2}
A Sturm permutation $\sigma\in S(N)$ is realizable as the Sturm permutation $\sigma=\sigma_g^\cN$ of some  spatially reversible nonlinearity $g\in\SO$, if and only if $\sigma$ is an integrable involution.
The same statement holds true in the Hamiltonian class $g\in\SH$.
In particular, any integrable Sturm involution $\sigma$ is realizable by a Hamiltonian nonlinearity $g=g(u)$.
\end{theorem}

\section{Neumann attractors and period maps}\label{Nlap}

\begin{figure}[t] 
\begin{center}
\includegraphics[width = \textwidth]{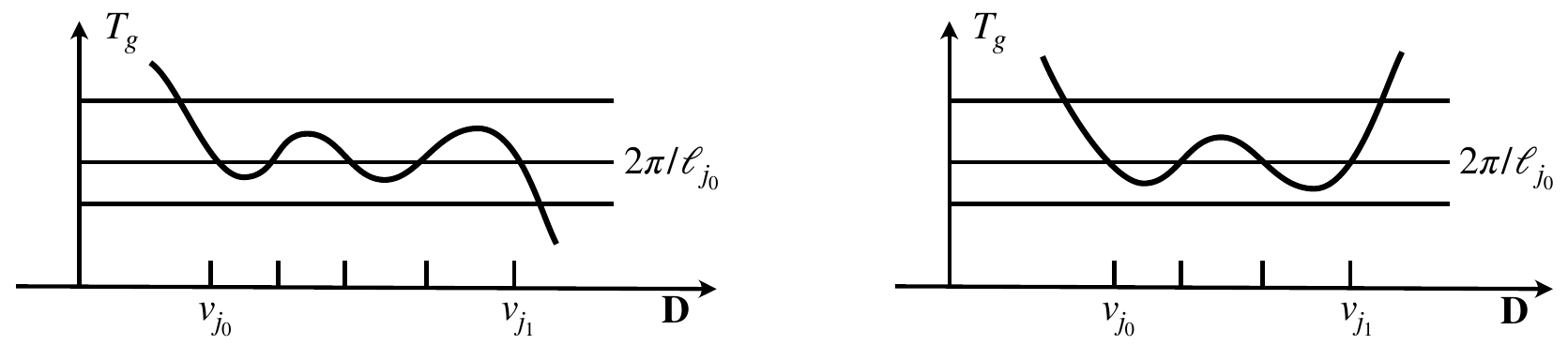}
\end{center}
\caption[Period maps]{\small\emph{ Examples of two period maps $T=T_g(v_0)$, partial view.
Horizontal axis: $v_0\in\mathbf{D}$\,. Vertical axis: $T$.
We assume that all values $T=2\pi/\ell,\ \ell\in\N$, are regular (horizontal lines),
i.e. the period maps $T$ cross such levels transversely, alternatingly up or down. 
In each connected component of $\mathbf{D}$, and at fixed $\ell$, each crossing indicates boundary values of the same type, say of local minimum type $v_0=\underline{v}_0$.
Note the odd number of consecutive crossings at locally transient mid-level $\ell$ (left), and the even number of crossings at locally extremal $\ell$ (right).}}
\label{figTg}
\end{figure}

In this section we investigate and compare the connection graphs $\cCN$ and $\cC_g^\cP$ of spatially reversible $g\in\SO$ via their ODE period maps, alias time maps, $T_g$ of \eqref{cT}.
Definition \ref{deflapsig} and theorem \ref{th3} formalize and characterize the relevant properties of the period map, via the abstract notion of a full lap signature $\mathfrak{S}=\mathfrak{S}(T_g)$.
The full lap signature is an expanded, but equivalent, version of the lap signature introduced in \cite{firowo12a}.
Lemma \ref{lemsigsig} establishes the one-to-one correspondence between full lap signatures and integrable Sturm involutions.
Theorem \ref{th4} therefore concludes that the full lap signature $\mathfrak{S}(T_g)$ also determines the (distinct, but closely related) connection graphs $\cCN$ and $\cC_g^\cP$ under Neumann and periodic boundary conditions.

\sss{Period maps}\label{T}

For spatially reversible nonlinearities $g\in\SO$, the minimal periods $T_g$ of periodic ODE orbits $(v,p)$ become an alternative to the characterization of global Neumann attractors $\cAN$ by integrable Sturm involutions $\sigma_g$\,.
We restrict the period map $T_g$ of \eqref{cT} to $\mathbf{D}\subseteq\R$, the intersection of the cyclicity set $\mathbf{C}$ with the horizontal $v$-axis:
\begin{equation}
\label{Tg}
T_g:\quad \mathbf{D}\rightarrow\R_+ \,.
\end{equation}
From \eqref{cT} we inherit the boundary of $\mathbf{D}$:
either $T_g\nearrow\infty$ blows up, or else the boundary point is an ODE center $v_0=e$ with $g_v(e,0)>0$ and finite limiting period 
\begin{equation}
\label{Te}
T_g\rightarrow2\pi\big/\sqrt{g_v(e,0)}\in(0,\infty)\,.
\end{equation}

The frozen waves $\mathbb{S}^1v\in\cF_g^\cP$ on the interval $0\leq x\leq 2\pi$, we recall, are in one-to-one correspondence \eqref{vpair} with the Neumann frozen min/max pairs $\underline{v}_0,\,\overline{v}_0$ which solve 
\begin{equation}
\label{Tv0}
T_g(v_0)=2\pi/\ell(v_0)\,.
\end{equation}
See figure \ref{figTg} for illustration.

The hyperbolicity condition \eqref{HYP} on frozen equilibria becomes 
\begin{equation}
\label{hypT}
T'_g(v_0)\neq0\,.
\end{equation}
In particular, hyperbolicity \eqref{HYP} holds if all  values $2\pi/\ell,\ \ell\in\N$, of $T_g$ are regular values.
This is a generic condition on $g$.
Indeed, perturbations $g^\eps(v,p):=(1+\eps)^{-2}g\bigl(v,(1+\eps)p\bigr)$ for small $|\eps|$ achieve $T_{g^\eps}=(1+\eps)T_g$\,.
Together with $(g_v(e,0))^2\not\in\N\cup\{0\}$ at all homogeneous equilibria $e\in\cE$, this establishes genericity of hyperbolicity \eqref{HYP} in the spatially reversible Neumann cases $g\in\SO$ and $g\in\SH$.

For frozen Neumann equilibria on $(0,\pi)$, we can refine hyperbolicity \eqref{hypT} to imply
\begin{equation} \label{iN}
i^\cN(\underline{v}_0) = 
\begin{cases} 
\,\ell\,, & \ \mbox{ if } \ T'_g(\underline{v}_0)<0\,; \\ 
\,\ell+1\,, & \ \mbox{ if } \ T'_g(\underline{v}_0)>0 \,.
\end{cases}
\end{equation}
Under periodic boundary conditions on $(0,2\pi)$, the Morse indices of frozen waves $v$ are 
\begin{equation} \label{iP}
i^\cP(v) = 
\begin{cases} 
\,2\ell-1 & \ \mbox{ if } \ T'_g(\underline{v}_0)<0\,; \\
\,2\ell & \ \mbox{ if } \ T'_g(\underline{v}_0)>0 \,.
\end{cases}
\end{equation}
See \cite{firowo04, firowo12a, firowo12b} for details. 
From now on we suppress dependence on the spatially reversible nonlinearity $g\in \SO$.

\sss{Full lap signatures}\label{S}

An essential tool for the classification of global attractors $\cA^\cN$ and connection graphs $\cC^\cN$ is the \emph{lap signature of the period map} $T=T(v_0)$. 
We proceed along the lines of section 5 in \cite{firowo12a} in equivalent, but somewhat expanded notation.

The min/max pairing $\underline{v}, \overline{v}$ of ODE periodic orbits decomposes the real domain $\mathbf{D}$ of $T$ in \eqref{Tg} into connected components
\begin{equation} \label{Dd}
\mathbf{D} = \bigcup_{1\le k\le K} (\,\mathbf{\underline{D}}_k \cup \mathbf{\overline{D}}_k\,),
\end{equation}
i.e.~into pairs of maximal intervals.
By the implicit function theorem, the ODE periodic orbits starting in the same interval are nested.
The $q$ corresponding min/max pairs $\underline{v}_0,\,\overline{v}_0$ define the min/max pairings $\mathbf{\underline{D}}_k, \mathbf{\overline{D}}_k$\,.
We call a min/max pair \emph{annular}, if $T\nearrow\infty$ at all four boundary points.
By dissipativity, the only other option are \emph{central} min/max pairs, where the interval pair shares an ODE center $e$, at a finite positive limit of $T$.
In the cyclicity set, central min/max pairs correspond to punctured disks.

In any $\mathbf{\underline{D}}_k$ let the list
\begin{equation}
\label{Sk}
\underline{S}_k:=\{\ell^k_1,\ldots,\ell^k_{s_k}\}
\end{equation}
collect the lap numbers $\ell(v_0)$ of \eqref{Tv0}, sorted by ascending order of $v_0$\,.
With the analogous definition for $\overline{S}_k$ in $\mathbf{\overline{D}}_k$\,, nesting implies the reverse order of the same lap numbers.
By construction, the list sizes $s_k$ add up to the total number $q$ of frozen waves, i.e. 
\begin{equation}
\label{qsk}
q=\sum_k s_k\,.
\end{equation}

To emphasize the min/max pairings $\underline{S}_k,\,\overline{S}_k$ we employ parenthesis symbols ``('' and ``)'' as follows.
In order of $v_0$\,, we replace each list $\underline{S}_k$ by the symbol ``(''.
The symbol ``('' reminds us, mnemonically, of  the local shape of any ODE periodic orbit $(v,p)$ through $(\underline{v}_0,0)\in\mathbf{\underline{D}}_k$ in the $(v,p)$-plane, as it crosses the horizontal $v$-axis at the minimum $\underline{v}_0$ of $v$.
Analogously, ``)'' indicates how local maxima are affiliated to $\overline{S}_k$, in reverse, throughout the paired upper interval $\mathbf{\overline{D}}_k$\,.
Only to central min/max pairs $\mathbf{\underline{D}}_k, \mathbf{\overline{D}}_k$ do we affiliate \emph{innermost parentheses} ``()'' even when $\underline{S}_k,\,\overline{S}_k$ happen to be empty.
For all other, i.e. annular, min/max pairs we require nonempty lists $\underline{S}_k,\,\overline{S}_k$ to grace them with parentheses.

This defines a \emph{regular parenthesis structure} of balanced open-close pairs $\ldots(\ldots)\ldots$, alias a \emph{Dyck word}, as follows.
Any annular open-close pair of parentheses, on the $v$-axis, comes with an associated upper arc of an ODE periodic orbit joining them, in the $(v,p)$-plane.
Thus it is the nesting of periodic orbits which defines the nesting and balancing, for pairs of parentheses.
Innermost central pairs () are trivial.
Note that redundant double parentheses $\ldots((\ldots))\ldots$ have to be forbidden, because the associated annulus or annuli would have to be nonempty, defining two successive connected components $\mathbf{\underline{D}}_k\,,\  \mathbf{\underline{D}}_{k+1}$\,.
Then $T\nearrow\infty$ at annulus boundaries requires at least one ODE saddle and one ODE center to separate the two parenthesis pairs.

To account for all intersections $v_0$ of the shooting meander $\cM$ in \eqref{shoot} with the $v$-axis, we actually augment the parenthesis and lap list structure as follows.
We insert bullet symbols $\bullet$\,, to account for any $\sigma$-stable PDE equilibria, i.e. for any ODE saddles $e\in\cE^\cN$.
Due to dissipativity, we prepend and append a leading and trailing symbol $\bullet$\,. 
These account for the minimal and maximal homogeneous equilibria $\underline{e},\overline{e}$, alias $\sigma(1)=1$ and $\sigma(N)=N$.
To avoid overprints, we move each min/max pair $\underline{S}_k,\,\overline{S}_k$ just inside the parenthesis pair which it defines:
\begin{equation}
\label{parS}
\mathfrak{S}(T)\ :=\ \bullet\,(\ldots (\underline{S}_k \ldots \overline{S}_k)\ldots )\,\bullet \,.
\end{equation}
Each central min/max pair $\underline{S}_k,\,\overline{S}_k$, possibly empty, we separate by a center symbol $\odot$ to account for the homogeneous PDE-unstable ODE center $e$ with $f_v(e,0)>0$ at the puncture between the central domain pair $\mathbf{\underline{D}}_k, \mathbf{\overline{D}}_k$.
All non-extremal ODE saddles appear exactly as configurations $\ldots )\,\bullet\,(\ldots$\,, alternatingly to innermost central min/max pairs $\ldots(\underline{S}_k \odot\overline{S}_k)\ldots$\,.
We call the resulting symbol sequence $\mathfrak{S}(T)$ the \emph{full lap signature} of the period map $T$.

The min/max pairings $\underline{v}_j\,, \overline{v}_j$ in corresponding order-reversed lists $\underline{S}_k\,, \overline{S}_k$ of the full lap signatures $\mathfrak{S}(T)$ keep track of the identification relation $\sim$ required in the application of corollary \ref{corCPN}.
The same lists of lap numbers also determine the signs of the period map derivatives $T'(v_0)$ required to determine the Morse indices of frozen equilibria by \eqref{iP}.
Indeed, we start from $T'(\underline{v}_1)<0$ in any list $\underline{S}_k$\,.

Our definition differs slightly from, but is equivalent to, a previous minimalistic notation in \cite{firowo12a, firowo12b, roc25, rofi14} which we called \emph{lap signature}.
The main difference is the ``redundant'' duplication of lists $\underline{S}_k$ by their reversals $\overline{S}_k$\,, and the explicit insertion of ``obvious'' locations $\bullet,\odot$ for spatially homogeneous ODE saddles and centers $e\in\cE$.

The following example with $n=5,\ q=4$ arises in figure \ref{figcyc}:
\begin{equation} \label{304}
\mathfrak{S}(T)\ =\ \bullet\Bigl(\{1,1\}\Bigl(\{1\}\odot\{1\}\Bigr)\,\bullet\Bigl(\{1\}\odot\{1\}\Bigr)\{1,1\}\Bigr)\,\bullet \,.
\end{equation}
Indeed note the $N=n+2q=13$ non-parenthesis symbolic entries $\bullet,\odot,1$ which indicate the horizontal order of boundary values at $x=0$ for homogeneous $\sigma$-stable ODE saddles $\bullet$\,, homogeneous PDE-unstable ODE centers $\odot$\,, and lap numbers.
The lap numbers, here $\ell=1$, appear at min/max pairs \eqref{vpair} with minimal ODE period $2\pi/\ell=2\pi$, and are located symmetrically to their $\sigma$-stable core, here the middle ODE saddle $\bullet$\,.
The annular domains $\mathbf{\underline{D}}_1$, $\mathbf{\overline{D}}_1$ of the outermost pair $(\ldots)$ are affiliated to the outer cyclicity region $\mathbf{C}_1$ in figure \ref{figcyc}
The central domains $\mathbf{\underline{D}}_2, \mathbf{\overline{D}}_2$ and $\mathbf{\underline{D}}_3, \mathbf{\overline{D}}_3$ of the two nonempty central pairs $\bigl(\{1\}\odot\{1\}\bigr)$ feature $\mathbf{C}_2$ and $\mathbf{C}_3$, respectively.

\sss{Abstract signatures, realization, and connection graphs}\label{SC}

We now formalize the notion of a full lap signature $\mathfrak{S}$, independently of period maps.
\begin{definition}\label{deflapsig}
Let $\mathfrak{S}$ be a finite sequence of symbols $\bullet,\odot,(,)$ and a collection of pairs $\underline{S}_k,\,\overline{S}_k$ of finite ordered lists $S=\{\ell_1,\ldots,\ell_{s_k}\}$ of positive integers $\ell_j\in\N$.
We call $\mathfrak{S}$ a \emph{full lap signature}, if the following nine properties (i)--(ix) all hold.
We first collect the five properties of symbols $\bullet,\odot,(,)$.
\begin{enumerate}[(i)]
	\item The sequence starts and ends with the symbol $\bullet$\,.
	\item The parenthesis symbols \emph{``(''} and \emph{``)''}, by themselves, form a regular parenthesis structure of balanced open-close pairs.
 		In particular, redundant parentheses $\ldots ((\ldots ))\ldots$ do not occur.
 		We call innermost pairs ``$()$'' \emph{central} and noncentral pairs \emph{annular}.
	\item Adjacent parentheses \emph{``)(''}  take the complete form $\ldots)\,\bullet\,(\ldots$ .
	\item Each central pair \emph{``()''} takes the complete form $\ldots(\underline{S}_k \odot\overline{S}_k)\ldots$ 			with possibly empty lists $\underline{S}_k\,,\,\overline{S}_k$\,.
	\item Each annular pair takes the complete form $\ldots(\underline{S}_k (\ldots)\,\bullet\ldots\bullet\,(\ldots)\overline{S}_k)\ldots$ with
		non\-empty lists $\underline{S}_k\,,\,\overline{S}_k$\,.
\end{enumerate}
Any pair of nonempty lists $S=\underline{S}_k=\{\ell_1,\ldots,\ell_{s_k}\}$ and $\overline{S}_k$ satisfies all of the following remaining four conditions.
\begin{enumerate}[(vi)]
	\item Each list $\overline{S}_k$ is the paired list $\underline{S}_k$ in reverse.
	\item \emph{Boundary conditions:}
		\begin{equation} \label{307}
		\ell_1=1\quad \textrm{and, in case of annular pairs: } \quad \ell_{s_k}=1 \textrm{ and } s_k \textrm{ is even.}
		\end{equation}
	\item \emph{Neighbor jump condition:}
		\begin{equation} \label{308}
		|\ell_{j+1}-\ell_j| \le 1 \,,\quad \textrm{for}\quad j=1,\dots,s_k-1 \,. \end{equation}
	\item[(ix)] \emph{Alternate jump condition:} \\ 
		If $\ell_{j_0-1} \ne \ell_{j_0} = \dots = \ell_{j_1} \ne \ell_{j_1+1}$ for $1 \le j_0 \le j_1<s_k$, then
		\begin{equation} \label{309}
		(\ell_{j_0}-\ell_{j_0-1})(\ell_{j_1+1}-\ell_{j_1})=(-1)^{j_1-j_0} \,.
		\end{equation}
		Here we define $\ell_0:=0$ in case $j_0=1$. 
\end{enumerate} 
\end{definition}

Continuity of the period map $T$ in \eqref{Tg}, and transversality \eqref{Tv0}, \eqref{hypT} to the regular values $T=2\pi/\ell$ impose strong restrictions on the lap signature of $T$. 
These restrictions in fact imply that the full lap signature $\mathfrak{S}(T)$ of $T=T_g$ is indeed a full lap signature $\mathfrak{S}$, in the sense of this abstract definition, for any $g\in \SO$.
See figure \ref{figTg} again for an illustration, in particular of alternate jump condition (ix).
The following theorem actually states how, conversely, the lap signatures $\mathfrak{S}(T)$ of such period maps cover all abstract lap signatures $\mathfrak{S}$.
For a proof see proposition 3 in \cite{roc07}.

\begin{theorem} \label{th3}
The full lap signature $\mathfrak{S}(T)$ of \,$T=T_g$ as defined in \eqref{Tg}--\eqref{parS} is a full lap signature $\mathfrak{S}$ in the sense of definition \ref{deflapsig}, for any $g\in \SO$.\\
Conversely, let $\mathfrak{S}$ denote any full lap signature in the abstract sense of definition \ref{deflapsig}. Then there exists $g\in \SO$ such that $\mathfrak{S}(T_g)=\mathfrak{S}$.
\end{theorem}

Regarding the connection equivalence of global attractors $\cA_f$\,, the main result of \cite{roc25} 
is the following.

\begin{theorem} \label{th4}
Let $g,\tilde{g}\in\SO$ possess the same full lap signature $\mathfrak{S}(T_g)=\mathfrak{S}(T_{\tilde{g}})$. 
Then their connection graphs are isomorphic:
\begin{equation} \label{312}
\cC_g^\cN\cong\cC_{\tilde{g}}^\cN, \quad \textrm{and} \quad \cC_g^\cP\cong\cC_{\tilde{g}}^\cP \,. 
\end{equation} 
\end{theorem}

\sss{Bijection between signatures and integrable Sturm involutions}\label{SigSig}

For the convenience of the reader, we present a proof of theorem \ref{th4} based on the results collected so far.
See also \cite{firo96, firo00, roc07, roc25}. 
The proof consists of a one-to-one translation between full lap signatures $\mathfrak{S}=\mathfrak{S}(T_g)=\mathfrak{S}(T_{\tilde{g}})$ and integrable Sturm involutions $\sigma_g=\sigma_{\tilde{g}}$\,.

\begin{lemma}\label{lemsigsig}
There is a bijective correspondence
\begin{equation}
\label{sigsig}
\mathfrak{S} \ \leftrightarrow\ \sigma
\end{equation}
between full lap signatures $\mathfrak{S}$ and integrable Sturm involutions $\sigma$.
\end{lemma}

\begin{proof}
We start from a full lap signature $\mathfrak{S}$ and construct an integrable Sturm involution $\sigma$.
We only have to identify the 2-cycles $(\underline{c}_{\,\alpha}\ \overline{c}_\alpha)$ of $\sigma$; see definition \ref{defint}.
By theorem \ref{th3}, there exists $g\in \SO$ such that the period map $T=T_g$ realizes $\mathfrak{S}=\mathfrak{S}(T)$.
Let $\sigma_g$ denote the associated integrable Sturm involution \eqref{202}, which is simply given by the permutations \eqref{lparity}.
In other words, all boundary values are fixed under $\sigma_g$\,, except for the 2-cycles which switch min/max pairs $\underline{v}_j\,, \overline{v}_j$ of odd lap numbers.

For example, let us enumerate the $N=9$ boundary values for the full lap signature $\mathfrak{S}$ which will arise as case 5.2(iv) of table \ref{tab21}:
\begin{equation}
\label{parse52iv}
\begin{array}{c c c c c c c c c c c c c c c c c}
      \mathfrak{S}&=&\bullet & \big(\,\{ & 1 &,& 1 & \}\big( & \odot & \big) \bullet \big( & \odot & \big)\{ & 1 &,& 1 & \}\,\big) & \bullet    \\[2mm]
      \mathfrak{S}^*&=&\bullet &  & 1 & & 1 & & \odot & \bullet & \odot & & 1 & & 1 &  & \bullet    \\[2mm]
         && 1 &  & 2 & & 3 &  & 4 & 5 & 6 &  & 7 & & 8 &  & 9  
\end{array}
\end{equation}
Since all lap numbers $\ell=1$ are odd,  it is exactly the two pairs $\underline{v}_j\,, \overline{v}_j= 2, 8$ and $3, 7$ which generate a 2-cycle, each.
This identifies the associated integrable Sturm involution as
\begin{equation}
\label{sigma52iv}
\sigma_g=(2\, 8)\,(3\, 7)= \{1,8,7,4,5,6,3,2,9\}\in S(9)\,.
\end{equation}

Abstractly, let $\mathfrak{S}^*$ denote the stripped symbol sequence $\mathfrak{S}$ which omits all parentheses ``('', ``)'', all list braces ``$\{$'', ``$\}$'', and all commas.
The stripped sequence of length $N=n+2q$ contains $n$ symbols $\bullet$ and $\odot$\,, for ODE saddles and ODE centers $(e,0)$, alias homogeneous PDE equilibria $e\in\cE_g$\,.
The remaining $2q$ entries are the Neumann lap numbers $\ell$ on $0<x<\pi$ associated to the $q$ spatially nonhomogeneous frozen PDE min/max pairs $(\underline{v}_j,\,\overline{v}_j)$ in $\cFN$, alias the  extrema of the associated $q$ frozen waves $v\in\cF_g^\cP$.
Consider the permutation $\sigma\in S(N)$ which swaps their ordinals, in the enumeration $\mathrm{id}=\{1,\ldots,N\}$ of $\mathfrak{S}^*$, for precisely those min/max pairs  $\underline{v}_j\,, \overline{v}_j$ which correspond to odd lap number entries of $\mathfrak{S}^*$.
Then $\sigma$ is an integrable  integrable Sturm involution.

Conversely, by theorem \ref{th2}, any integrable Sturm involution $\sigma$ is realized as $\sigma=\sigma_g$\,, for some nonlinearity $g\in \SO$.
Let $T=T_g$ denote the associated period map.
Then we have already explicated how to arrive at the associated full lap signature $\mathfrak{S}:=\mathfrak{S}(T)$.
See theorem \ref{th3}.
Alternatively, see section 6 of \cite{firowo12a} for the completely explicit construction of the lap signature from any integrable Sturm involution.

Since the two constructions are mutually inverse, this proves the lemma.
\end{proof}

\emph{Proof of theorem} \ref{th4}.\quad
By lemma \ref{lemsigsig}, equal lap signatures $\mathfrak{S}(T_g)=\mathfrak{S}(T_{\tilde{g}})$ imply equal Sturm permutations $\sigma_g=\sigma_{\tilde{g}}$\,.
By theorem \ref{thfCAN}(ii), this implies $\cCN\cong\cC_{\tilde{g}}^\cN$\,.
Since $\mathfrak{S}$ also keeps track of min/max pairings $\sim$, uniquely, theorem \ref{thgCPN} implies $\cC_g^\cP\cong\cC_{\tilde{g}}^\cP$\,; see \eqref{gCPN}.
This proves theorem \ref{th4}.
$\hfill\bowtie$

In view of characterization theorem \ref{th3} and realization theorem \ref{th2}, the global attractors, Sturm meanders, and signature classes can all be realized by Hamiltonian period maps of $g\in\SH$, in the Neumann case.
The proof of theorem \ref{th4} in \cite{roc25}, via direct homotopies of period maps, nevertheless requires homotopies in the larger spatially reversible class $g\in\SO$, so far.
We also caution the reader once again that, even in the Hamiltonian or spatially reversible cases, the orbit equivalence of the global attractors $\cA_g^\cP$ and $\cA_{\tilde{g}}^\cP$, analogously to \eqref{312}, has \emph{not} been proved.

\section{Proof of theorem \ref{th21}}\label{PfConn}

\begin{table}[t!]
\begin{center}
\begin{tabular}{ l  c c c c } \hline
$n.q$ & $n+q$ & $\#\odot$ & $N=n+2q$ & full lap signature $\mathfrak{S}$ \\ \hline \\
1.0     & 1     & 0  &    1           & $\bullet$ \\ \\ 
3.0     & 3     & 1  &    3           & $\bullet\,\big(\odot\big)\,\bullet$ \\  \\ 
3.1     & 4     & 1  &    5           & $\bullet\,\big(\,\{1\}\odot\{1\}\,\big)\,\bullet$ \\ \\ 
3.2(i)  & 5     & 1  &    7           & $\bullet\,\big(\,\{1,2\}\odot\{2,1\}\,\big)\,\bullet$ \\ 
3.2(ii)      &       &      &                 &  $\bullet\,\big(\,\{1,1\}\odot\{1,1\}\,\big)\,\bullet$ \\ \\ 
3.3(i)  & 6     & 1  &    9           &  $\bullet\,\big(\,\{1,2,3\}\odot\{3,2,1\}\,\big)\,\bullet$ \\ 
3.3(ii)      &       &      &                 &  $\bullet\,\big(\,\{1,2,2\}\odot\{2,2,1\}\,\big)\,\bullet$ \\ 
3.3(iii)     &       &      &                 &  $\bullet\,\big(\,\{1,1,1\}\odot\{1,1,1\}\,\big)\,\bullet$ \\ \\               
3.4(i)  & 7     & 1  &   11          & $\bullet\,\big(\,\{1,2,3,4\}\odot\{4,3,2,1\}\,\big)\,\bullet$ \\
3.4(ii)      &       &      &                 & $\bullet\,\big(\,\{1,2,3,3\}\odot\{3,3,2,1\}\,\big)\,\bullet$ \\
3.4(iii)     &       &      &                 & $\bullet\,\big(\,\{1,2,2,2\}\odot\{2,2,2,1\}\,\big)\,\bullet$ \\
3.4(iv)     &       &      &                 & $\bullet\,\big(\,\{1,2,2,1\}\odot\{1,2,2,1\}\,\big)\,\bullet$ \\                  
3.4(v)      &       &      &                 & $\bullet\,\big(\,\{1,1,1,2\}\odot\{2,1,1,1\}\,\big)\,\bullet$ \\ 
3.4(vi)     &       &      &                 & $\bullet\,\big(\,\{1,1,1,1\}\odot\{1,1,1,1\}\,\big)\,\bullet$ \\ \\
5.0     & 5     & 2  &   5            & $\bullet\,\big(\odot\big)\bullet\big(\odot\big)\,\bullet$ \\ \\
5.1     & 6     & 2  &   7            & $\bullet\,\big(\,\{1\}\odot\{1\}\,\big)\bullet\big(\odot\big)\,\bullet$ \\ \\
5.2(i)     & 7     & 2  &   9         &$\bullet\,\big(\,\{1,2\}\odot\{2,1\}\,\big)\bullet\big(\odot\big)\,\bullet$ \\
5.2(ii)          &       &      &             & $\bullet\,\big(\,\{1,1\}\odot\{1,1\}\,\big)\bullet\big(\odot\big)\,\bullet$ \\
5.2(iii)         &       &      &             & $\bullet\,\big(\,\{1\}\odot\{1\}\,\big)\bullet\big(\,\{1\}\odot\{1\}\,\big)\,\bullet$ \\                  
5.2(iv)        &       &      &              & $\bullet\,\big(\,\{1,1\}\big(\odot\big) \bullet \big(\odot\big)\{1,1\}\,\big)\,\bullet$  \\ \\
7.0     & 7     & 3  &   7            & $\bullet\,\big(\odot\big)\,\bullet\,\big(\odot\big)\,\bullet\,\big(\odot\big)\,\bullet$ \\ \\ 
\hline
\end{tabular} 
\end{center} 
\caption[Full lap signatures for $n+q\leq7$]{The last column enumerates all full lap signatures $\mathfrak{S}$ of period maps $T$ with $n+q\leq7$, up to trivial equivalences. 
Here $n$ is the odd number of homogeneous equilibria, $\#\odot=(n-1)/2$ is the number of ODE centers $\odot$\,, and $q$ is the number 
of frozen waves after the application of the freezing and symmetrizing homotopies.
The list is grouped by case labels $n.q$, lexicographically.
The total number $N$ of Neumann equilibria, frozen or homogeneous, is $n+2q$, because each frozen wave $v$ on $x\in\bS^1=\R/2\pi\Z$ generates a pair $\underline{v},\, \overline{v}$ of Neumann equilibria on the half interval $0<x<\pi$; see \eqref{vpair}.}
\label{tab21}
\end{table}

In this section we prove theorem \ref{th21} on the 21 connection graphs $\cCP$ of the Sturm global attractors $\cAP$ of parabolic PDEs \eqref{101}, under the restriction $n+q\leq 7$.
See figure \ref{figgraphs}.
We recall that $n$ counts the homogeneous equilibria $e\in\cEP$ and $q$ the frozen or rotating waves $v\in\cFP\cup\cRP$.
The nonlinearities $f$ are restrained to the $\bSO$ Sturm class $\SSOP$, i.e. $f=f(v,p)$ is $C^2$, dissipative, and all homogeneous equilibria, frozen and rotating waves are hyperbolic; see \eqref{HYP}.

By theorem \ref{thgCPN} and corollary \ref{corCPN}, the result is equivalent to determining the Neumann connection graphs $\cCN/\!\!\sim$ in the further restricted $\bO$ class of spatially reversible nonlinearities $g\in\SO$, under the same restriction $n+q\leq7$.
The Neumann connection graphs $\cCN$ are always determined by their Sturm permutations $\sigma_g$; see theorem \ref{th1}.
In the spatially reversible class of $g$, the Sturm permutations are integrable involutions $\sigma_g$; see definition \ref{defint} and theorem \ref{th2}.
By Lemma \ref{lemsigsig}, integrable involutions $\sigma$ stand in bijective correspondence to full lap signatures $\mathfrak{S}$, abstractly.
Theorems \ref{th3} and \ref{th4} realize this correspondence, concretely, as a correspondence between Sturm permutations $\sigma_g$ and full lap signatures $\mathfrak{S}=\mathfrak{S}(T_g)$ of period maps $T_g$\,, in the class of spatially reversible $g\in\SO$.

It is therefore sufficient to enumerate the abstract full lap signatures $\mathfrak{S}$ with $n$ symbols $\bullet,\odot$\,, corresponding to the $n$ homogeneous equilibria, and with min/max parenthesized paired lists $\underline{S}_k\,,\,\overline{S}_k$ totaling $q$ lap number entries $\ell$ in the lower lists $\underline{S}_k$\,; see definition \ref{deflapsig}.
We proceed along table \ref{tab21} to obtain the 21 resulting signatures $\mathfrak{S}$ with $n+q\leq7$, in the last column.

\subsection{Enumeration of full lap signatures}\label{Tab21}

The first column $n.q$ of table \ref{tab21} lists the cases by increasing odd $n=1,3,5,7$; see section \ref{Intro}.
The trivial case $n=1$ consists of a single homogeneous equilibrium $\underline{e}=\overline{e}=e$ as the global attractor $\cA=\{e\}$, which does not accommodate any parentheses in $\mathfrak{S}=\bullet$\,. 
In particular $q=0$.
For each $n>1$, we group the cases by increasing $q=0,\ldots,7-n$.
Roman numerals indicate subcases.
The second and fourth columns list the sum $n+q$, and the total number $N=n+2q$ of non-parenthesis entries in $\mathfrak{S}$, respectively.
This also determines the total number $N$ of Neumann boundary values $e,\underline{v}_j\,, \overline{v}_j$ which appear in the associated integrable Sturm involution $\sigma\in S(N)$; compare \eqref{201}, \eqref{202}.

Saddles $\bullet$ and centers $\odot$ alternate with each other. 
Indeed, properties (iii) and (iv) in definition \ref{deflapsig} imply that all parentheses have to open, after any $\bullet$ and up to the next $\odot$\,.
Analogously, parentheses have to close, after any $\odot$ and up to the next $\bullet$\,.
By property (i), this leaves us with the counts $\#\odot=(n-1)/2$ in column three.

To arrive at the remaining 20 nontrivial full lap signatures $\mathfrak{S}$ of the last column, let us address $n=3$ first.
Then $\#\odot=1$, and the two saddles $\bullet$ must be leading and terminating $\mathfrak{S}$.
Since property (ii) forbids redundantly repeated parentheses, this only accommodates one single pair of parentheses.
Therefore, any $\mathfrak{S}$ must take the form 
\begin{equation}
\label{n=3}
\bullet\big(\underline{S}\odot\overline{S}\big)\bullet
\end{equation}
with a single order-reversed min/max pair $\underline{S},\overline{S}$ of lap lists.
See properties (ii), (iv), (vi) of definition \ref{deflapsig}, and note that the min/max pair is central.
Also note that each of the reversely ordered lists $\underline{S}$ has length $q$.
For each $0\leq q\leq4$, we have listed the options compliant with properties (vi)--(ix) of $\underline{S}$, in inverse lexicographic order.
This settles all 13 cases with $n=3$.

The cases $n=5$ possess $\#\odot=2$ centers.
The centers alternate with three saddles $\bullet$: one leading, one terminating, and the third separating the two centers.
This provides the general template
\begin{equation}
\label{n=5}
\bullet\big(\underline{S}_1\big(\underline{S}_2\odot\overline{S}_2\big)\bullet\big(\underline{S}_3\odot\overline{S}_3\big)\overline{S}_1\big)\bullet
\end{equation}
with up to one nonempty annular min/max list pair $\underline{S}_1,\overline{S}_1$, and two (potentially empty) central min/max list pairs associated to 
$\underline{S}_2$ and $\underline{S}_3$\,, respectively.
Compare figure \ref{figcyc} and \eqref{304}, where all three lists are nonempty, albeit at the expense of $q=4$ and slightly too large $n+q=9$.

The presence of a nonempty annular min/max pair $\underline{S}_1,\ \overline{S}_a$ (and its parentheses) requires $0<s_1\leq q\leq 2$ to be even; see definition \ref{deflapsig}(v),(vii).
Hence $s_1=2$, and boundary condition \eqref{307} implies $\underline{S}_1=\{1,1\}$.
Since \eqref{qsk} does not leave room for any other $s_k$\,, the central lists $\underline{S}_2$ and $\underline{S}_3$ are both empty.
This provides case 5.2(iv) of table \ref{tab21}.

Now all remaining cases of $n=5$ can only feature central min/max pairs, possibly empty.
The general template \eqref{n=5} for them reduces to
\begin{equation}
\label{n=5c}
\bullet\big(\underline{S}_2\odot\overline{S}_2\big)\bullet\big(\underline{S}_3\odot\overline{S}_3\big)\bullet\ .
\end{equation}
For $q=0$, the lists $\underline{S}_2\,, \underline{S}_3$ must both be empty; see case 5.0.
For $q=1$, one list is empty and, by boundary condition (vi), the other one-element list must coincide with $\underline{S}=\{1\}$; see case 5.1.
We have omitted the reflected case 
\begin{equation}
\label{5.1'}
\bullet\,\big(\odot\big)\bullet\big(\{1\}\odot\{1\}\big)\,\bullet\ .
\end{equation}
Indeed, any realization of \eqref{5.1'} leads back to case 5.1, by flow equivalence under the linear substitution $u\mapsto -u$, alias $g(v,p)\mapsto -g(-v,-p)$.
On the meander level, we have called this a \emph{trivial equivalence} of global attractors and connection graphs; see 
(6.7)--(6.9) in \cite{firo96}.
By construction, this trivial equivalence reverses the full lap signature $\mathfrak{S}$ -- which we omit, for brevity.

\begin{figure}[t!] 
\begin{center}
\includegraphics[width = \textwidth]{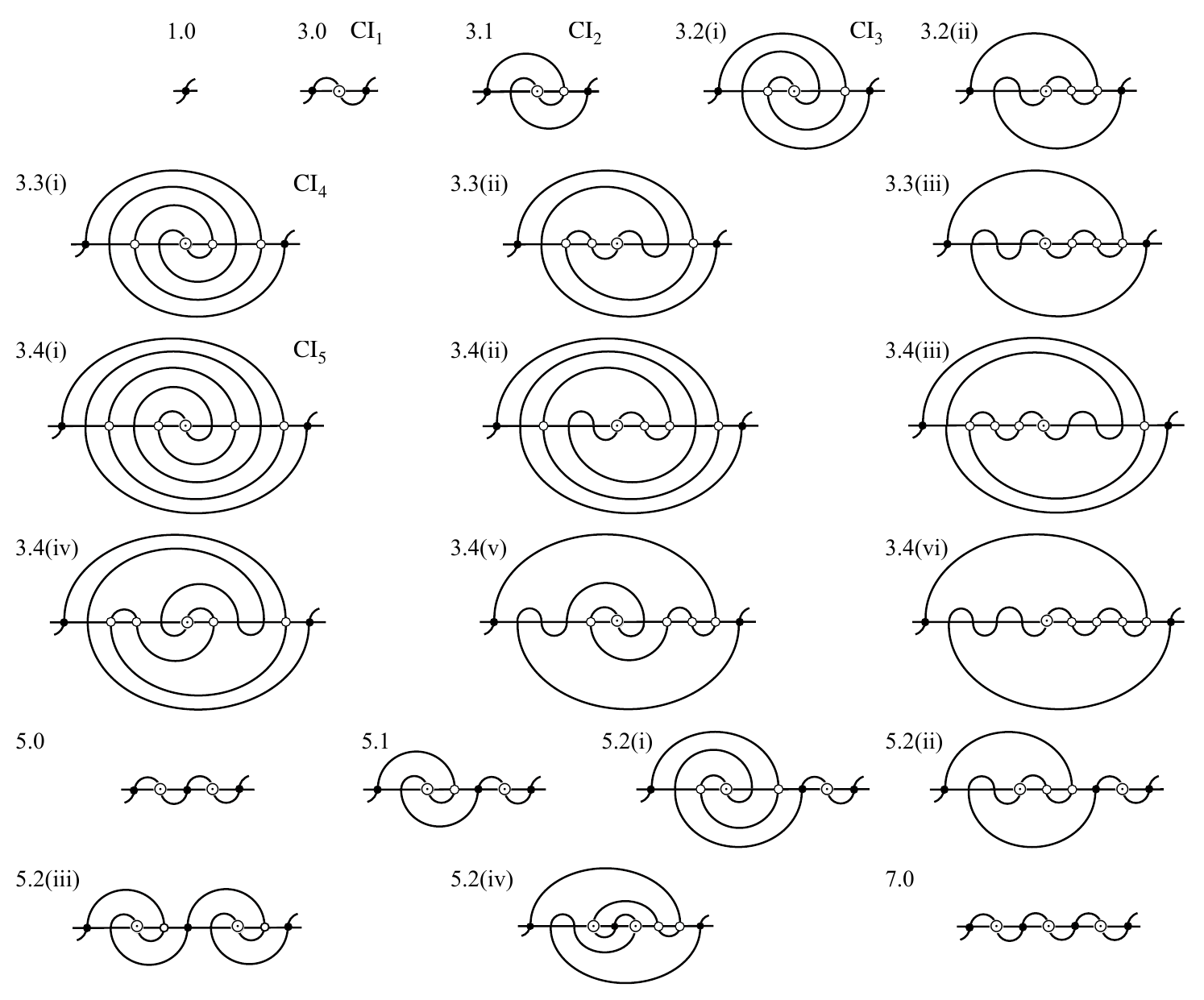}
\end{center}
\caption[Neumann meanders for $n+q\leq7$]{\small\emph{ List of all 21 stylized Neumann meanders $\cM_f$ corresponding to the full lap signatures $\mathfrak{S}$ of table \ref{tab21}. 
The meanders are based on the derivation of the associated integrable Sturm involutions $\sigma$, as in \eqref{parse52iv}, \eqref{sigma52iv}.
The listing is lexicographic in $n.q$, for $n+q\leq7$, with $n$ homogeneous equilibria and $q$ pairs of frozen equilibria. 
See the legends of figure \ref{figcyc} and table \ref{tab21} for symbols $\bullet,\,\odot,\,\circ$, and figure \ref{figgraphs} for the resulting connection graphs.
}}
\label{figmeanders}
\end{figure}

The non-annular cases 5.2 of $s_2+s_3=q=2$ offer two options.
First suppose $s_2=2,\ s_3=0$, up to trivial equivalence.
Properties (vii)--(ix) of definition \ref{deflapsig} only admit $\underline{S}_2=\{1,2\}$ and $\underline{S}_2=\{1,1\}$.
This establishes cases 5.2(i) and 5.2(ii), respectively.
The final option $s_2=1=s_3$ provides the last remaining case 5.2(iii).

The final case $n=7$ is trivial because $q=0$ does not leave room for any nonempty lap lists.
This proves table \ref{tab21}. \hfill $\bowtie$

\subsection{From lap signatures to Sturm permutations}\label{Sig2Sig}

The bijection \eqref{sigsig} of lemma \ref{lemsigsig} translates any full lap signature $\mathfrak{S}(T_g)$ to the associated integrable Sturm involution $\sigma_g$.
In \eqref{parse52iv}, \eqref{sigma52iv} we have already worked out the details of this translation for case 5.2(iv) of table \ref{tab21}.
All remaining full lap signatures $\mathfrak{S}=\mathfrak{S}(T_g)$ of table \ref{tab21} parse to their associated integrable Sturm involutions $\sigma_g$\,, by the same recipe.
Figure \ref{figmeanders} represents the resulting permutations by stylized meander curves $\cM_g$\,.
The involutions $\sigma_g$ appear along the horizontal axis, if we enumerate the intersection points along the meander; see \eqref{shoot}.
Reverting the meander orientation, by $180^\circ$ rotation $(v,p)\mapsto(-v,-p)$, produces the omitted trivially equivalent Sturm permutations and Neumann attractors.

\subsection{From Sturm permutations to connection graphs}\label{Sig2C}

\begin{figure}[t] 
\begin{center}
\includegraphics[width = 0.35\textwidth]{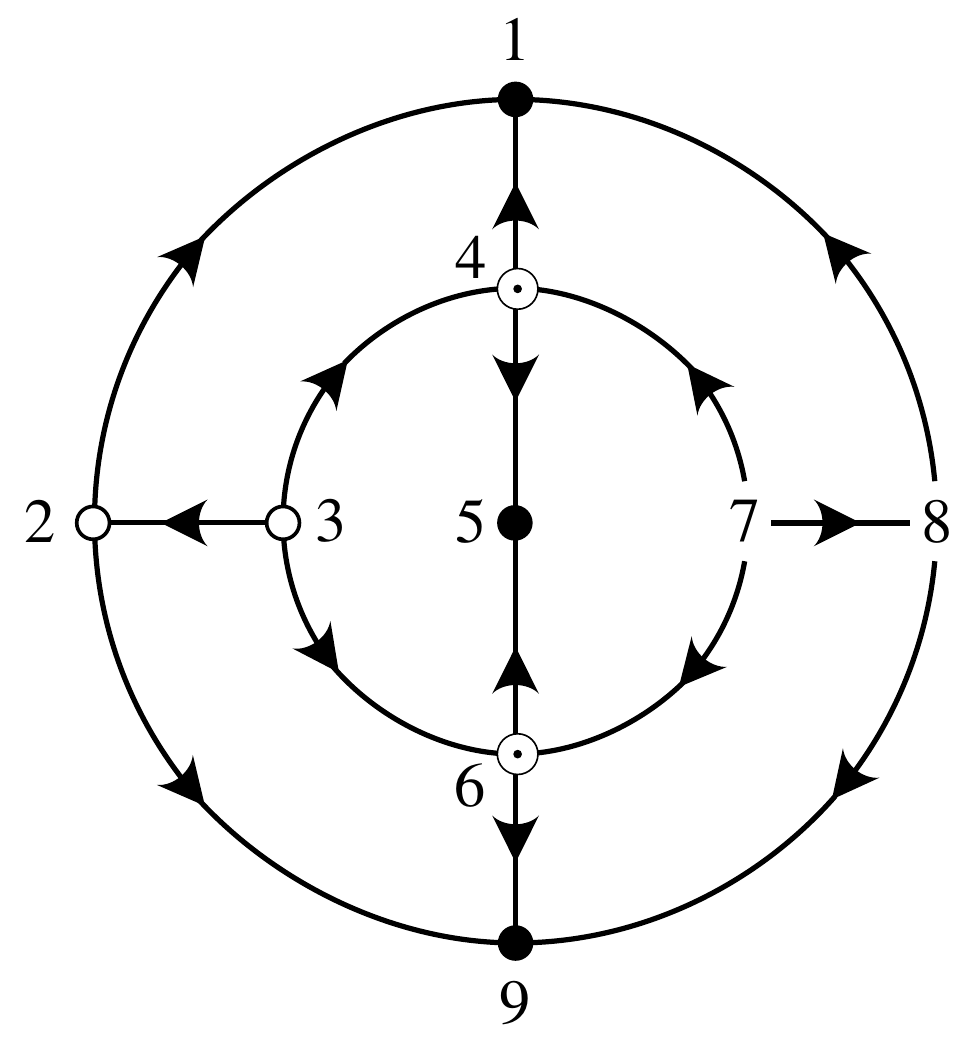}
\end{center}
\caption[Neumann connection graph of example 5.2(iv)]{\small\emph{ The Neumann connection graph $\cCN$ associated to the integrable Sturm involution $\sigma$ of case 5.2(iv) in table \ref{tab21} and figure \ref{figmeanders}.
Numbering and symbols follow the order along the meander in figure \ref{figmeanders}.
The PDE-stable ODE saddles $1,5,9\ (\bullet)$ and the PDE-unstable ODE centers $4,6\ (\odot)$ are spatially homogeneous, $n=5$.
The $q=2$ frozen wave pairs $\underline{v}_j\,,\overline{v}_j$ are $2,8$ and $3,7$, with minima $\underline{v}_j$ marked by circles.
The diagram first appeared in \cite{fiedlertatra}. 
Rotation around the vertical axis gives an impression of an associated three-dimensional global attractor $\cAP$ with two rotating waves, one through each frozen pair.
}}
\label{fig52iv}
\end{figure}

The passage from Sturm permutations $\sigma_g$ to Neumann connection graphs $\cCN$ has been established almost 30 years ago; see Theorem \ref{thfCAN}(ii) and \cite{firo96}.
By \eqref{ij}, the Sturm permutation determines all Morse indices $i$ and all zero numbers $z$ of differences among the stationary PDE solutions.
The connection graph then follows by $\cN$-adjacency and blocking \eqref{block}.
See \cite{fiedlertatra} for a complete classification of all Sturm permutations and connection graphs with $N\leq9$ stationary solutions.
Our case 5.2(iv), for example, figures there as $\sigma=\pi_{17}$ in case 2.4 and figure 3.5; see our figure \ref{fig52iv}.
See also \cite{firo10} for a classification of all planar Neumann cases with $N=11$ stationary solutions, and \cite{firo18} for a classification of all 31 Sturm 3-balls with $N\leq13$.

We invoke corollary \ref{corCPN} to pass from the Neumann connection graphs $\cCN$ to the connection graphs $\cC_g^\cP\cong\cCN/\!\!\sim$ under periodic boundary conditions.
In our example 5.2(iv) we have to identify the min/max vertex pairs $\underline{v}_j\,, \overline{v}_j= 2, 8$ and $3, 7$ of the Neumann connection graph in figure \ref{fig52iv}.
Obviously, this leads to case 5.2(iv) in figure \ref{figgraphs}, where the two circles denote the two identified min/max paired vertices, respectively.
Their unstable dimensions follow from the associated period map $T_g$ in the annular region; see \eqref{iP} with $\ell=1$.
The unstable dimensions of homogeneous equilibria are easily determined explicitly, from sines and cosines.
A tempting template for the geometry of the global attractor $\cA_g^\cP$ appears by rotation of figure \ref{fig52iv} around the vertical axis of spatially homogeneous solutions.
Checkably, this provides the correct connection graph.
And we have issued sufficient warnings against jumping to conclusions for $C^0$ orbit equivalence of global $\bSO$-attractors, by now.

The remaining meanders of figure \ref{figmeanders} can be treated analogously, to arrive at \ref{figgraphs}.
This completes the proof of theorem \ref{th21}.\hfill $\bowtie$

\begin{figure}[t] 
\begin{center}
\includegraphics[width = 0.75\textwidth]{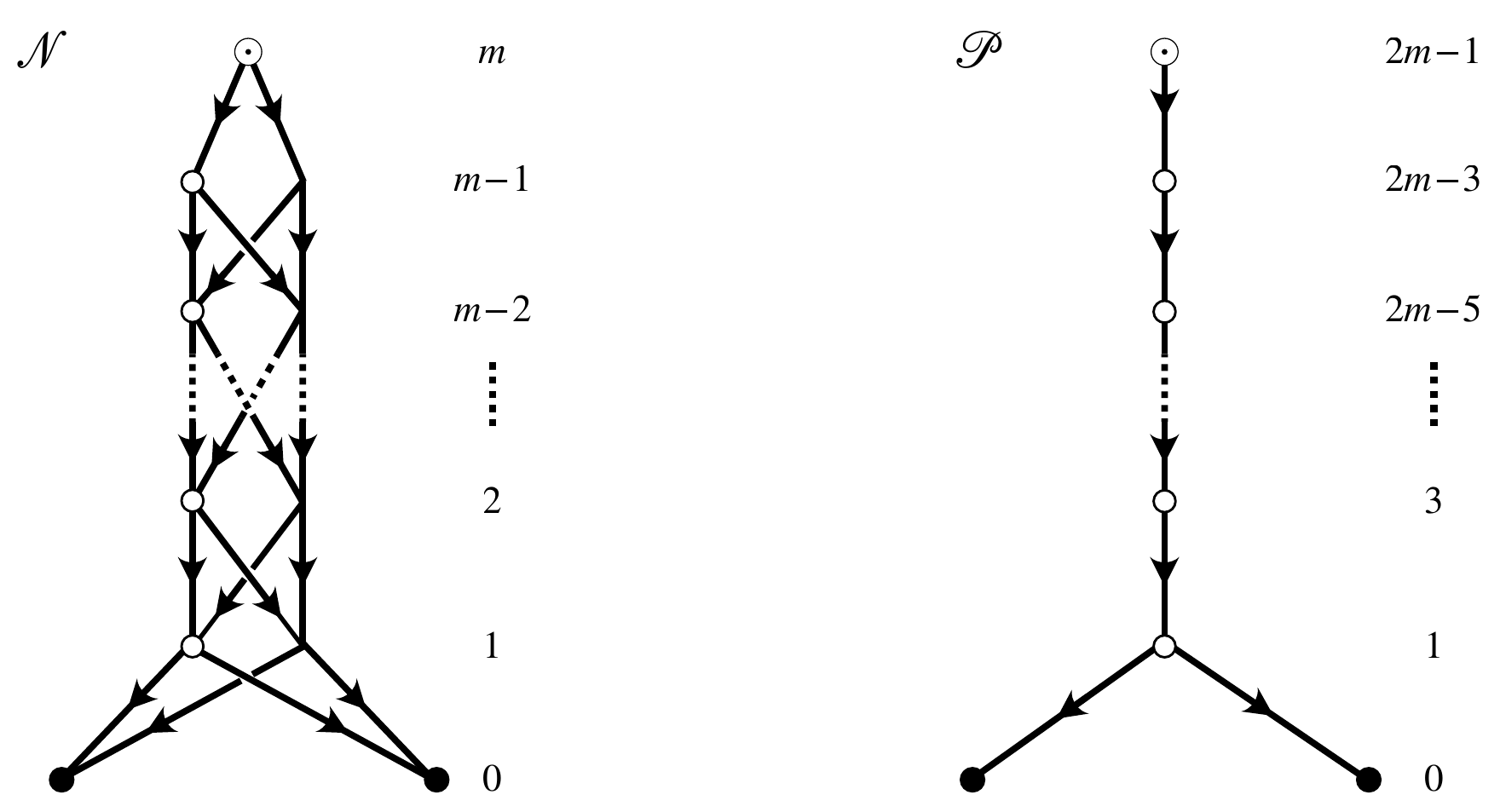}
\end{center}
\caption[Connection graphs of Chafee-Infante attractors]{\small\emph{ The connection graphs  $\cC(\mathrm{CI}_m^{\cN,\cP})$ of the Chafee-Infante global attractors $\mathrm{CI}_m^{\cN,\cP}$.
Left: Neumann boundary $\cN$, dimension $m$.
Right: Periodic boundary conditions $\cP$, dimension $2m-1$.
Note the $n=3$ homogeneous equilibria:
two of them PDE-stable ODE saddles $\bullet$, at $e=\pm1$, and the single ODE center $\odot$, at $e=0$, of maximal Morse indices $i=m$ and $i=2m-1$, respectively.
The periodic case $\cP$ (right) has $q=m-1$ frozen waves $v_j=\circ,\ j=1,\ldots,m-1$, of (half period) lap numbers $\ell=j$ and Morse indices $i=2j-1$. See \eqref{iP}.
In the Neumann case $\cN$ (left), each $v_j$ is represented by a min/max pair $\underline{v}_j=\circ,\ \overline{v}_j$ (not marked), of equal lap number and Morse index $\ell=j=i$.
See \eqref{vpair}, \eqref{iN}.
We may just view $\mathrm{CI}_m^\cP$ as a rotation figure of $\mathrm{CI}_m^\cN$\,.
}}
\label{figCI}
\end{figure}

\subsection{Example: Chafee-Infante attractors}\label{secCI}
As examples in any dimension $m\geq1$, we mention the celebrated Chafee-Infante attractors $\mathrm{CI}_m^\cN=\cAN$; see  \cite{chin74, hen81}.
They arise as global attractors of the Neumann PDE \eqref{PDEg} for Hamiltonian nonlinearities $g=g(v)=\lambda\,v(1-v^2)$ and $(m-1)^2<\lambda<m^2$.
See \cite{firo22} and the references there for a detailed discussion. 
Their full lap signatures are 
\begin{equation}
\label{CImlapsign}
\mathfrak{S}(\mathrm{CI}_m^\cN)\ =\ \bullet\,\big(\{1,\ldots,m-1\}\odot\{m-1,\ldots,1\}\big)\,\bullet\,.
\end{equation}
The $q=m-1$ frozen equilibrium pairs identify this as the lap signatures of cases 3.$(m$-$1)$ or 3.$(m$-$1)$(i) in figures \ref{figgraphs}, \ref{figmeanders}, and table \ref{tab21}, for Chafee-Infante dimensions $m=1,2,3,4,5$.
Note the ODE saddles $\bullet$ at $e=\pm1$, and the ODE center $\odot$ of PDE Morse index $i=m$ at $e=0$.
The symmetric pairs of frozen equilibria are nested around $v=e=0,\ p=0$ in the phase plane \eqref{vpg}.
They possess lap numbers and Morse indices $\ell=i=\{1,\ldots,m-1\}$, from outside inward.
The resulting meander, or simple symmetry arguments, show that any two equilibria of adjacent Morse levels connect heteroclinically.
The equivalence $\sim$ identifies the two frozen equilibria at the same Morse level.
Under periodic boundary conditions on $0\leq x\leq 2\pi$, corollary \ref{corCPN} therefore implies the Chafee-Infante connection graph $\cC(\mathrm{CI}_m)^\cP$ of figure \ref{figCI}.
For unstable dimensions compare \eqref{iN}, \eqref{iP} with $T'<0$.

\subsection{Pitchforkability}\label{Pitch}
Sturm permutations $\sigma$ were introduced by Fusco and Rocha in \cite{furo91}.
The same paper explores how to track Sturm permutations and connection graphs through pitchfork bifurcations, for Neumann boundaries.
For the original idea see \cite{hen81,hen85}.

Sturm permutations $\sigma\in S(N)$ are called \emph{pitchforkable} if the total number $N\geq3$ of stationary solutions can be reduced to $N-2$ by a (reverse) pitchfork bifurcation.
Formally, this requires three consecutive entries to be adjacent, ascending or descending, among the sequence $\sigma=\{\sigma(1)\ \ldots\ \sigma(N)\}$.
Alas, it turns out that non-pitchforkable Sturm permutations do exist, and in fact abound.
The smallest examples require $N=11{:}\ \sigma=(1\ 6\ 7\ 10\ 3\ 4\ 9\ 8\ 5\ 2\ 11)$ and $\sigma=(1\ 8\ 7\ 2\ 3\ 6\ 9\ 10\ 5\ 4\ 11)$; see \cite{firo96}.

Integrable Sturm involutions $\sigma$, in contrast, are always pitchforkable.
In fact they can be reduced to the trivial case $N=1$ by a sequence of formal pitchfork reductions.
Indeed, assume $N\geq3$.
Then the full lap signature $\mathfrak{S}$, associated to $\sigma$ by lemma \ref{lemsigsig}, contains at least one ODE center $\ldots\big(\underline{S}_k\odot\overline{S}_k\big)\ldots$\ .

In case $\underline{S}_k$ is empty, the signature continues with ODE saddles, i.e.~$\ldots\bullet(\odot)\bullet\ldots$, and $\sigma$ encounters three consecutively ascending fixed points.
Pitchfork reduction replaces this signature by $\ldots\bullet\ldots$, and the reduced sequence remains a lap signature; see definition \ref{deflapsig}.

In case $\underline{S}_k=\{\ell_1,\ldots,\ell_{s_k}\}$ is nonempty, the signature $\mathfrak{S}$ looks like $\ldots,\ell_{s_k}\}\odot\{\ell_{s_k},\ldots$\,, locally.
For even $\ell_{s_k}$\,, we encounter three consecutively ascending fixed points, as before.
For odd $\ell_{s_k}$\,, the fixed point $\odot$\,, say at position $m$, is framed by a 2-cycle $(m\!+\!1\ \ m\!-\!1)$, and we obtain three consecutively descending entries $\sigma=(\ldots\ m\!+\!1\ m\  m\!-\!1\ \ldots )$.
In either case, the reduced sequence remains a lap signature again.

On a formal level, this suggests a pitchfork approach to build up the connection graphs, and even the global attractors up to orbit equivalence, starting from scratch, i.e. from the trivial case $N=1$ upwards.
See \cite{roc25}.
Although complete technical details have not been pushed through, to our knowledge,
they look manageable for the spatially reversible Neumann case $g\in\SO$.
Alternatively, one might attempt to carry the pitchfork approach over to periodic boundary conditions $g\in\SO$ and aim for their connection graphs $\cC_g^\cP$, by direct induction.
More ambitiously, it is tempting to aim for full $C^0$ orbit equivalence of global attractors $\cA_g^\cP$ with the same full lap signature $\mathfrak{S}(T_g)$, in the class of spatially reversible $g$.

\section{Global dynamics in positive delayed feedback}\label{Delay}

Surprisingly, the parabolic global attractors presented so far may serve as a benchmark for the 
global dynamics of positive feedback delay differential equations (DDEs)
\begin{equation}\label{501}
\dot y(t) = h(y(t), y(t - 1)) \,.
\end{equation}
A sufficient dissipativity condition, for example, is $y_1\,h(y_1,y_2)<0$ for sufficiently large $|y_1|$.
We assume \emph{positive delayed feedback} $\partial_2 h > 0$ for $C^1$ dissipative $h{:}\ \R^2\to\R$.

DDEs \eqref{501} find wide applications, e.g.~in biology, to 
simulate delayed self-excitatory effects in ``population'' dynamics ranging from animal 
demographics to protein concentrations. 
The specific modeling context often differs fundamentally from the reaction-drift-diffusion processes behind PDE \eqref{101}.
Nevertheless, parabolic PDEs \eqref{101} and DDEs \eqref{501} exhibit remarkable similarities, 
as dynamical systems. 
We highlight three geometrical features.

First, like PDE \eqref{101}, the dissipative DDE \eqref{501} generates a compact solution semiflow on an infinite-dimensional state space $y(t+.)\in X:=C^0([-1,0],\R)$ of solutions \cite{halu93}. 
Hence, \eqref{501} possesses a compact global attractor $\cA_h^\cD$.
Unlike the parabolic case, though, even the linear delay semiflow is not analytic.

Second, the DDE \eqref{501} possesses a \emph{delay zero number}, similar to the parabolic one, and a Poincaré--Bendixson theorem ensues; see  \cite{mpse96}, but also contrast with \cite{furo24}.
Under hyperbolicity, the recurrent set of \eqref{501}  on $\cA_h^\cD$ consists of finitely many hyperbolic
equilibria $e\in\cE_h^\cD$, i.e. $h(e,e)=0$, and hyperbolic nonstationary periodic solutions $y=v\in\cR_h^\cD$;
compare \eqref{HYP}. 
Analogously to \eqref{cAP}, all transients are then heteroclinic connections $u\in\cH_h^\cD$ among them:
\begin{equation}
\label{cAD}
\cA_h^\cD=\cE_h^\cD\cup\cR_h^\cD\cup\cH_h^\cD\,.
\end{equation}
Third, the linearizations of \eqref{501} at $\cE_h^\cD$ and $\cR_h^\cD$ possess a 
Sturm-Liouville like eigenvalue structure which relates Morse indices and delay zero numbers of eigenfunctions \cite{mpnu13, lop23}.

\begin{figure}[t] 
\begin{center}
\includegraphics[width =0.5\textwidth]{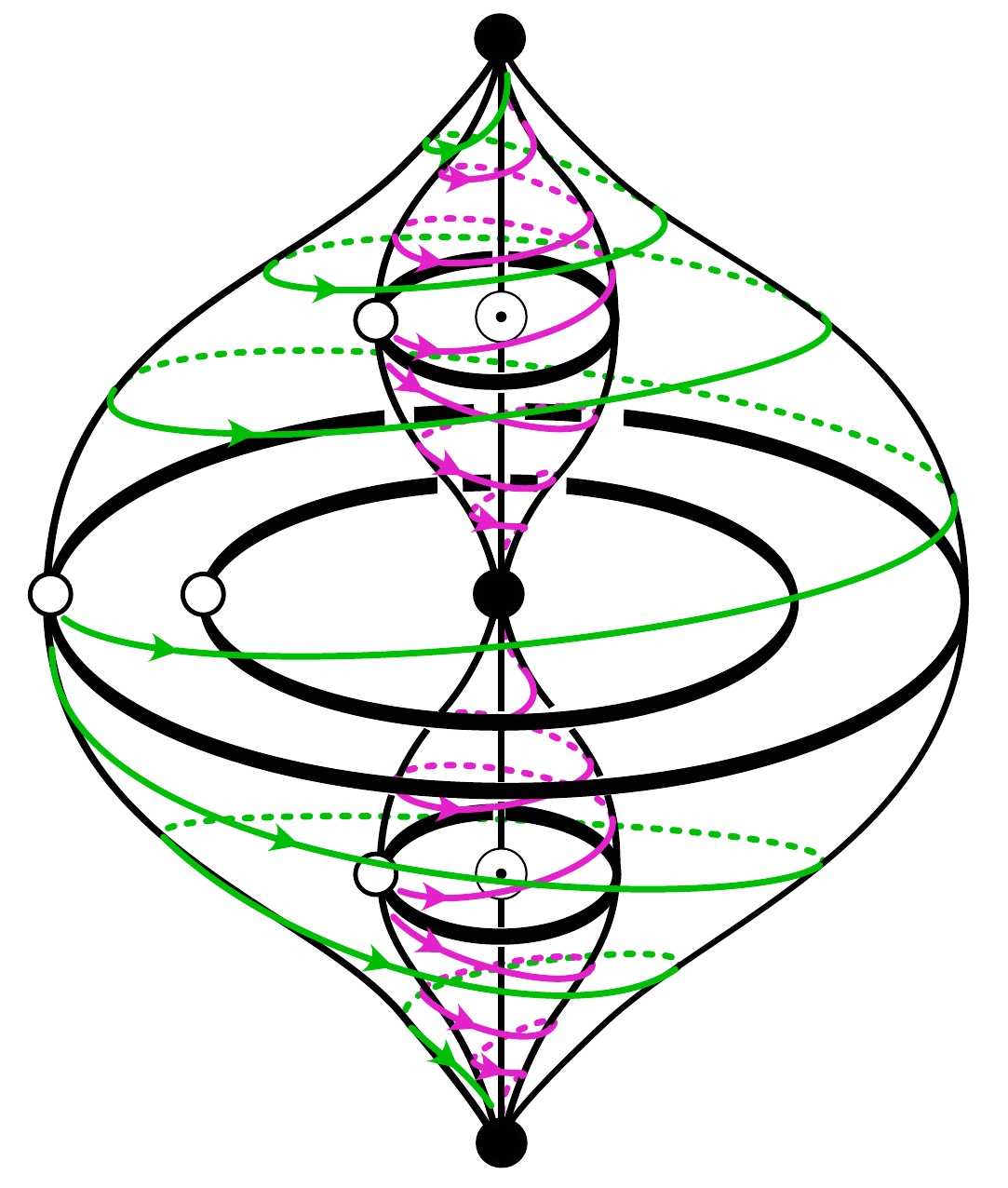}
\end{center}
\caption[The Vas tulip]{\small\emph{ The three-dimensional \emph{Vas tulip} $\AT$\,, a global attractor which arises in the context of dissipative DDEs \eqref{501} and was discovered by \cite{krvas11}.
It features $n=5$ equilibria $e\in\cE_\cT^\cD$ (stable: $\bullet$, unstable: $\odot$) and $q=4$ periodic orbits $v\in\cR_\cT^\cD\ (\circ)$.
Note the two three-dimensional Chafee-Infante spindles $\mathrm{CI}_2^\cP$, each with one periodic orbit of unstable dimension $i=1$, which are stacked on top of each other along the vertical ``rotation'' axis.
The two spindles are surrounded, at the waist, by two horizontally nested, large annular periodic orbits of Morse indices $i=1$ and $i=2$, respectively.
Green and purple: some heteroclinic orbits emanating from the three rotating waves with Morse index $i=1$ along their two-dimensional unstable ``separatrix'' manifolds. 
For the connection graph $\CT$ see also figure \ref{figNP}(right).
}}
\label{figtulip}
\end{figure}

These similarities embolden us to compare the \emph{delay connection graphs} $\cC_h^\cD$ of delay 
attractors $\cA_h^\cD$\,, for DDEs \eqref{501}, with the \emph{parabolic connection graphs} $\cCP$ of PDEs
\eqref{101}. 
In contrast to systematic parabolic studies of $\cCP$, we are aware of only three main example classes of $\cC_h^\cD$, in the literature on positive monotone feedback DDEs: 
\begin{enumerate}[Ex 1.]
	\item The \emph{spindle class} of $h$ in \cite{kr08,krwawu99} features $q=m-1$ periodic orbits nested around a single central equilibrium with Morse index $i=2m-1$, and two extremal stable equilibria $\underline{e},\,\overline{e}$. 
\mbox{Their delay} connection graph coincides with the Chafee-Infante connection graph $\cC(\mathrm{CI}_m^\cP)$ of dimension $2m-1$.
See figure \ref{figCI}(right).
	\item The three-dimensional \emph{Vas tulip} $\AT$ was the first DDE attractor with a connection graph outside the spindle class; see figure \ref{figtulip} and \cite{krvas11}. 
It contains $n=5$ equilibria and $q=4$ periodic orbits.
The Vas tulip consists of two three-dimensional Chafee-Infante spindles $\mathrm{CI}_2^\cP$\,, stacked on top of each other, and surrounded by two annular periodic solutions of large amplitude. 
The delay connection graph $\cC_\cT^\cD$ coincides with the parabolic connection graph $\mathcal{C}_g^\cP$ which we will design in figure \ref{figNP}.
	\item A third, rich source of examples with positive delayed feedback originates from DDE \eqref{501} with \emph{negative delayed feedback} $\partial_2h_- < 0$. 
	We then observe analogous zero number and spectral structures, but stable periodic orbits may arise. 
	The \emph{enharmonic} class in \cite{lop23} only restricts $h_-(v,p)$ in \eqref{501} to be even in $v$ and odd in $p$. As a consequence, the equilibria and periodic orbits of the DDE attractor $\cA_{h_-}^\cD$ can be fully understood in terms of a period map and its full lap signature, akin to $\mathfrak{S}$ in \eqref{Sk}, \eqref{parS} and definition \ref{deflapsig}. 
	In combination with an automatic transversality property analogous to \eqref{transv}, the delay attractor $\cA_{h_-}^\cD$ is Morse-Smale if all equilibria and periodic solutions are hyperbolic. 
	Hence, the enharmonic connection graphs arise by tracking sequences of Hopf bifurcations and obey blocking and adjacency rules which are quite analogous to the Neumann PDE attractors of section \ref{Nmeander}.

	Let us perturb a bounded, negative delayed feedback, enharmonic nonlinearity to 
	\begin{equation}\label{modifDDE} 
	h(v,p) = -\varepsilon v - h_-(v, p)\,.
	\end{equation}
	Then the resulting DDE \eqref{501} has positive delayed feedback $\partial_2 h > 0$ and is dissipative. 
	For $\varepsilon > 0$ sufficiently small, the delay connection graph $\cC^{\cD}_h$ turns out to be an unstable double-cone suspension of the negative delayed feedback connection graph. 
	The results in \cite{lop23} imply that the enharmonic modification \eqref{modifDDE} realizes the class $3.q$ of parabolic connection graphs with $n=3$ equilibria and $q$ periodic orbits; compare figure \ref{figgraphs} and table \ref{tab21}. 
	Notably, the Vas tulip above does not fit into the enharmonic class, because it requires $n=5$ equilibria. 
\end{enumerate} 
All above three classes of DDE examples with positive delayed feedback also occur as connection graphs $\cCP$ of the parabolic PDE \eqref{101}.
For the Vas tulip, see section \ref{PDET}.

Inspired by the parabolic PDE \eqref{101} and the Morse--Smale property of $\cA_h^\cD$, one strategy to construct further examples of delay connection graphs is to address local bifurcations. 
See our discussion of Neumann pitchforks, in section \ref{Pitch}.     
A major obstruction to the bifurcation approach is that the set $\cR_h^\cD$ of periodic orbits of the DDE \eqref{501} is not understood as well as in the parabolic case, globally. 
For example, \eqref{501} lacks any overt or hidden $\bSO$ or $\bO$ symmetries, which assisted our PDE analysis so much.
Nevertheless, Vas \cite{vas17} broke ground: she showed how periodic solutions of \eqref{501} realize certain parenthesis expressions that encode the relative position of delay periodic orbits. 
Her method is a first step towards characterizing all possible configurations of periodic orbits.
Unfortunately, it does not reveal any underlying bifurcation process towards them.  
The elegant Vas constructions excel, in fact, at identifying regions with periodic orbits, but do not count them there. 
Guided by the Vas method, \cite{lop25} has embarked on a synthesis with pitchforks of period maps, in the spirit of sections \ref{T} and \ref{Pitch}.
Via a second-oder ODE associated to delay periodic orbits of \eqref{501}, a combinatorial structure emerges, akin to the parabolic full lap signatures $\mathfrak{S}$ of definition \ref{deflapsig}.

The conditions (i)--(ix) of definition \ref{deflapsig} are minimal requirements for full lap signatures $\mathfrak{S}$ which are based on continuous period maps $T_g$.
Therefore, the parabolic PDEs \eqref{101} can be expected to replicate any sequence of local bifurcations in the DDEs \eqref{501}, under positive delayed feedback.
A forteriori, the parabolic classification approach of our present paper may become essential for the DDEs \eqref{501}.
In particular, we expect all delay connection graphs $\cC_h^\cD$ to be contained in the parabolic class $\cC_g^\cP$ of spatially reversible nonlinearities $g\in\SO$.

\section{Example: a Sturm realization of the frozen Vas tulip}\label{PDET}

In the previous section \ref{Delay} we have presented the Vas tulip attractor $\AT$ for the DDE \eqref{501}; see figure \ref{figtulip}.
In the present section, we set out to construct a frozen version of $\AT$ as a Sturm PDE global attractor $\cA_g^\cP$, in the spatially reversible class $g\in\SO$.
Since $g$ is spatially reversible, indeed, all $q=4$ periodic orbits of the Vas tulip attractor $\AT$ will have to be frozen in $\cA_g^\cP$.
In particular, the two global attractors cannot be orbit equivalent.
A posteriori, of course, we might invert the lines of section \ref{Freeze} and reanimate the four frozen waves into rotating waves, by an unfreezing homotopy.

Instead of such orbit equivalence, we just aim for a parabolic connection graph $\cC_g^\cP\cong\CT$ which coincides with the connection graph $\CT$ of the Vas tulip; see figure \ref{figNP}(right).
We will also show that $\CT$ corresponds to the lap signature \eqref{304}.

Brute force approaches to this realization problem might start from computer generated lists of Sturm permutations $\sigma\in S(N)$ for Neumann problems $g\in\SO$. 
For the Vas tulip with $N=n+2q=13$ we obtain 1083 Sturm candidates, to be scanned for the Vas tulip by the explicit methods of \cite{firo96} and theorem \ref{thfCAN}, algorithmically.
Alternatively, we could list all full lap signatures $\mathfrak{S}$ in the class $n.q=5.4$, analogously to table \ref{tab21}, and then test them all, as in sections \ref{Tab21}--\ref{Sig2C}.
Alas, such uninspired blunder contributes exactly nothing to our principal goal: geometric insight.

\sss{Heuristics}\label{Heu}

\begin{figure}[t] 
\begin{center}
\includegraphics[width =  \textwidth]{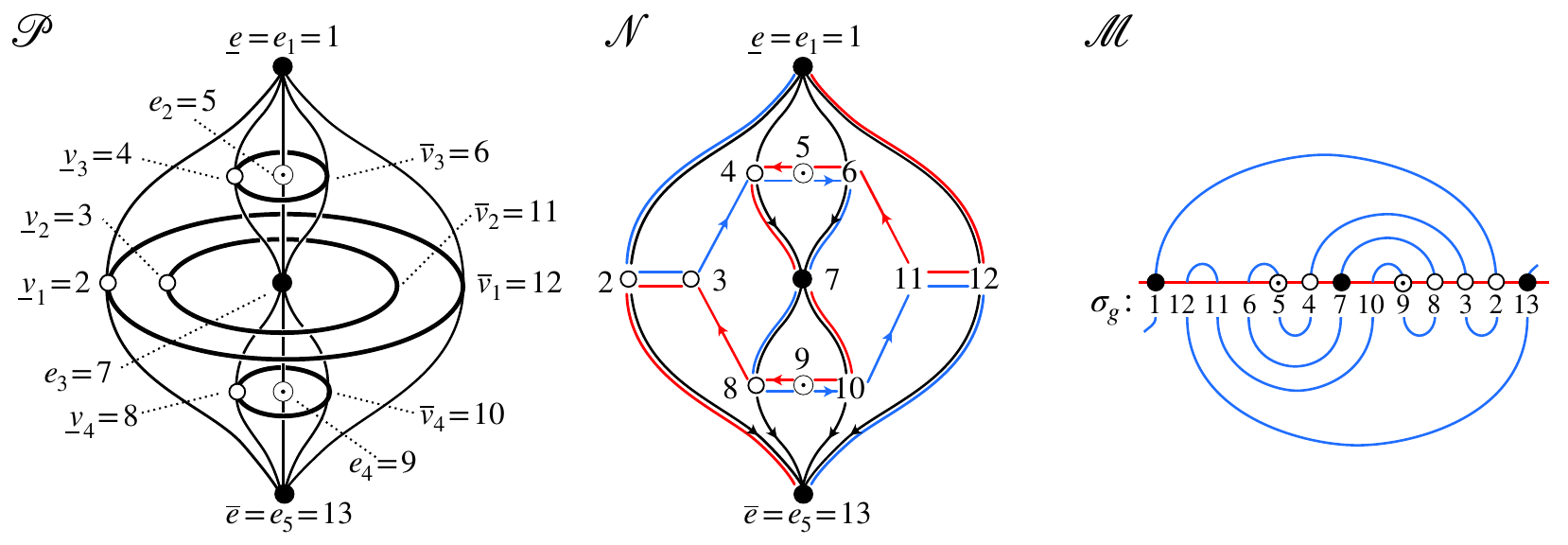}
\end{center}
\caption[Frozen Vas tulip, Neumann section, and meander]{\small\emph{ 
\emph{Left:}
The Vas tulip delay attractor $\AT$ of figure \ref{figtulip}, redrawn with labels which refer to the parabolic frozen tulip $\cA_g^\cP$.
Homogeneous equilibria: $\bullet$ (stable) and $\odot$ (unstable).
Frozen circles $v_j\in\cF_g^\cP$ replace the $q=4$ periodic orbits in $\cR_\cT^\cD$\,.
Circles $\circ$ mark their minima $\underline{v}_j$\,, left, opposite to their non-marked maxima $\overline{v}_j$\,.
\\
\emph{Center:}
Vertical planar section of the frozen Vas tulip, as a Thom-Smale complex, in same annota\-tion.
The six unstable separatrix manifolds of saddles $2,4,6,8,10,12$ with Morse index $i=1$ define the edges of a Thom-Smale dynamic complex with sink vertices $1,7,13\ (\bullet)$ and faces marked by $3,5,9,11$.
A \emph{bipolar orientation}, from \emph{north pole} $\underline{e}=1$ to \emph{south pole} $\overline{e}=13$, directs all edges downward (black arrows), non-dynamically.
The orientation induces Hamiltonian paths $h_0$ (blue) and $h_1$ (red), pole-to-pole; see text.
Then $\sigma_g:=h_0^{-1}\circ h_1\in S(13)$ is a Sturm permutation such that the connection graph $\cCN$ fits the prescribed Thom-Smale complex \cite{firo08,firo09,firo10,firo22}. 
Boundary orders \eqref{201}, \eqref{202} at $x=0$ and at $x=\pi$ are given by $h_0$ and $h_1$\,, respectively.
In our case, the two paths $h_0$ (blue) and $h_1$ (red) define the integrable Sturm involution $\sigma_g$ of \eqref{sigmaT}.
\\
\emph{Right:}
The stylized Sturm meander $\cM_g$ associated to the integrable Sturm involution $\sigma_g$.
Integer labels and symbols $\bullet,\circ,\odot$ as on left and center.
}}
\label{figTNM}
\end{figure}

Our construction of a  PDE version of the DDE Vas tulip $\AT$ will be geometric.
Intuitively, figure \ref{figtulip} looks like a \foreignlanguage{russian}{матрёшка} (matryoshka) of nested spindles, i.e.~a figure of rotation around the vertical axis -- which, in absence of $\bSO$ equivariance, it is not.
Still, we observe three DDE-stable equilibria (bullets $\bullet$) and two unstable ones (circled dots $\odot$) along the vertical ``rotation axis''.
Circles $\circ$ mark the four time-periodic orbits $y=v_j$ of delay equation \eqref{501} at their minima $\underline{v}_j$\,, on the left, opposite to their non-marked maxima.
See figure \ref{figTNM}(left) for a frozen version.

In figure \ref{figTNM}(center), we draw a schematic cross section of the \foreignlanguage{russian}{матрёшка}, with identical annotations.
Roughly speaking, the passage from the parabolic $\bSO$ connection graph $\cC_g^\cN$, under periodic boundary conditions, to the Neumann version $\cCN$ in corollary \ref{corCPN} looks just like a passage to such a cross section.
In the horizontal radial coordinate, the identification ``$\sim$'' would work like $r\sim-r$.
We now view the cross section diagram as a planar Thom-Smale complex, not unlike using radial components $r\in\R$ of polar coordinates for the horizontal plane in spite of lost again $\bSO$ equivariance.
The six unstable separatrix manifolds of saddles $2,4,6,8,10,12$ with Morse index $i=1$ define the edges of a Thom-Smale dynamic complex with sink vertices $1,7,13\ (\bullet)$ and faces marked by $3,5,9,11$.
We seek to realize  that Thom-Smale complex as the dynamics complex of unstable manifolds, associated to the planar Neumann attractor $\cAN$ of an integrable Sturm involution $\sigma_g$\,.

\sss{Planar Sturm design}\label{Des}

\begin{figure}[t] 
\begin{center}
\includegraphics[width = \textwidth]{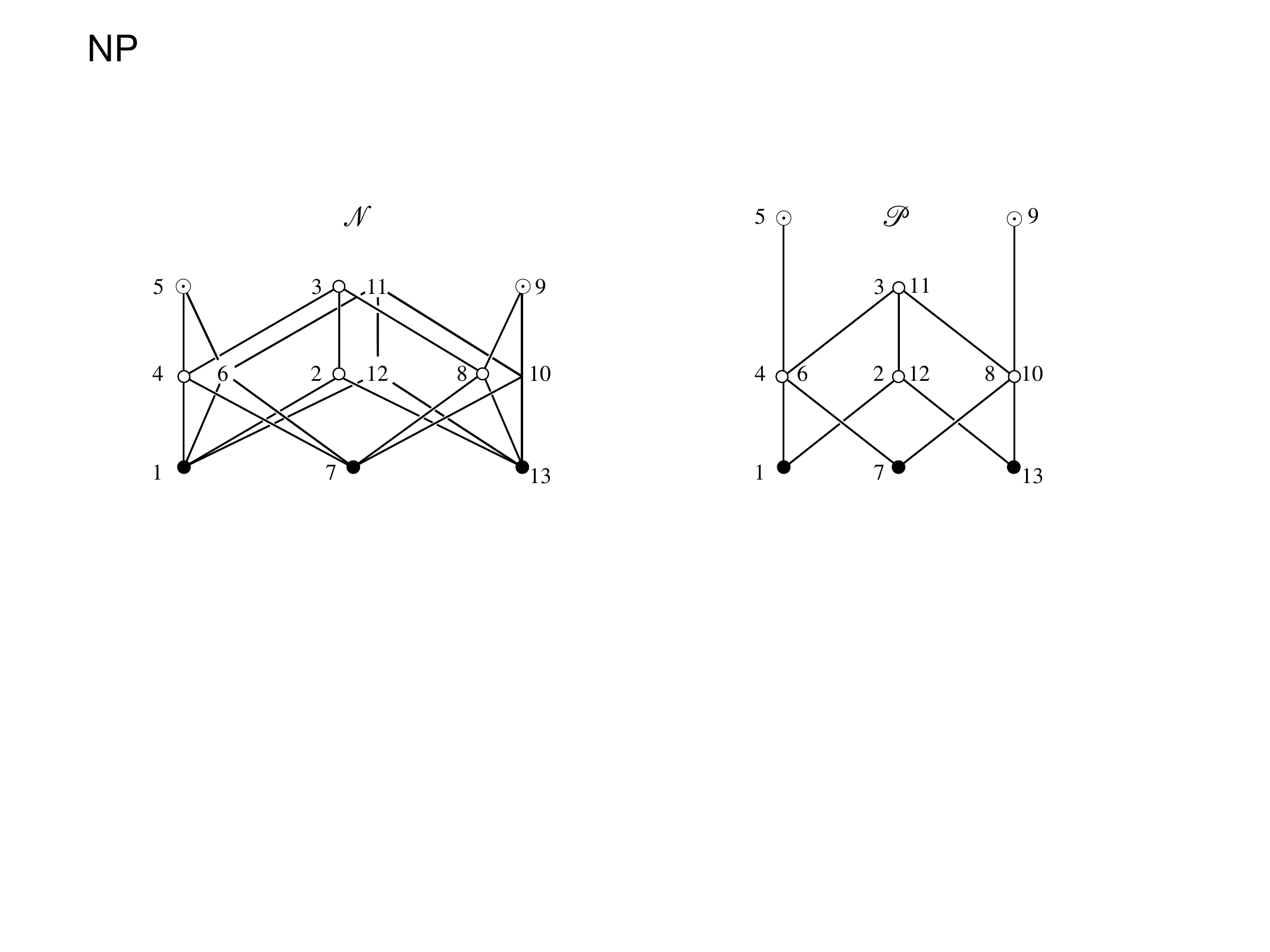}
\end{center}
\caption[The Vas tulip: PDE connection graphs]{\small\emph{ 
The Vas tulip: Neumann and periodic connection graphs.
For annotations and symbols see figure \ref{figTNM}.
\emph{Left:} 
The parabolic Neumann connection graph $\cCN$ of the meander $\cM_g$ of figure \ref{figTNM}(right).
Vertical grading is by Morse indices $i$.
Frozen equilibrium vertices are grouped into min/max pairs $\{2,12\},\ \{3,11\},\ \{4,6\},\ \{8,10\}$ by their $q=4$ ``hidden'' frozen waves; see \eqref{vpair}.
Each min/max pair corresponds to a 2-cycle of $\sigma_g$\,.\\
\emph{Right:}
The two vertices of each grouped frozen min/max pair have been identified to a single frozen wave $\circ$ by the relation $\sim$ of \eqref{CPN}.
This provides the parabolic connection graph $\cC_g^\cP$ of \eqref{101} with $f:=g$, under periodic boundary conditions.
By direct comparison with figures \ref{figtulip} and \ref{figTNM}(left), $\cC_g^\cP$ is indeed isomorphic to the delay connection graph $\CT$ of the Vas tulip delay attractor $\AT$ of delay equation \eqref{501}.
}}
\label{figNP}
\end{figure}

In \cite{firo10} we have presented many examples for the design of planar Sturm attractors $\cAN$ with a prescribed Thom-Smale complex, based on \cite{firo08, firo09}.
See also \cite{firo22}.
The main ingredient of design is an acyclic \emph{bipolar orientation} of the 1-skeleton, from the minimal \emph{north pole} equilibrium $\underline{e}=1$ to the maximal \emph{south pole} equilibrium $\overline{e}=N=13$: downward black arrows in figure \ref{figTNM}(center).
This orientation is \emph{not} dynamic, but corresponds to the monotone order $z=0$ along both branches of any one-dimensional unstable separatrix manifold.
For given bipolar orientation, we define pole-to-pole Hamiltonian paths $h_0$ (blue in figure \ref{figTNM}, center) and $h_1$ (red) from  $\underline{e}=1$ to $\overline{e}$ as follows.
Along separatrix edges, both paths follow the bipolar orientation.
The path $h_0$ traverses each 2-cell from the lowest left boundary saddle to the highest right boundary saddle: lower left to upper right.
The rule for the red path $h_1$, however, reads: lower right to upper left.
For example, $h_0{:}\  \ldots 2\ 3\ 4\ldots$ (blue) and $h_1{:}\  \ldots\ 8\ 3\ 2\ \ldots$ (red) through the same 2-disc face of vertex 3 with $i=2$. 
The main result of \cite{firo08, firo09} then asserts that $\sigma_g:=h_0^{-1}\circ h_1$ is a Sturm permutation such that
\emph{the Neumann connection graph} $\cCN$ \emph{coincides with the connection graph of the prescribed Thom-Smale complex.}
Since we have labeled the vertices such that the path $h_0=\mathrm{id}$ simply follows 1,\ldots,13, in ascending order, the red path $h_1$ in the center figure of \ref{figTNM} coincides with the realizing Sturm permutation:
\begin{equation}
\label{sigmaT}
\begin{aligned}
\sigma_g=h_0^{-1}\circ h_1  &=\{1,12,11,6,5,4,7,10,9,8,3,2,13\}\\
		&= (2\;12)\, (3\;11)\, (4\;6)\, (8\;10)\,.
\end{aligned}
\end{equation}
In figure \ref{figTNM}(right) we represent $\sigma_g$ by the associated stylized meander $\cM_g$\,.
The labeled intersections with the horizontal axis indicate boundary values at $x=\pi$ and correspond to vertices in the cross-sectional Thom-Smale complex and the frozen Vas tulip, with the same integer labels and symbols.
For the associated connection graph $\cC_g^{\cN}$ see figure \ref{figNP}(left), again with the same labels and symbols.

At first glance, the above approach only provides Neumann realizations $\sigma=\sigma_g$ by some $g\in\SxN$; see theorem \ref{thfCAN} and \eqref{sturmNx}.
Direct inspection of \eqref{sigmaT}, however, reveals that $\sigma_g$ is an integrable Sturm involution; see definition \ref{defint}.
This provides realizations $\sigma=\sigma_g$ by spatially reversible $g\in\SO$, and even by Hamiltonian nonlinearities $g\in\SH$; see theorem \ref{th2}.
The integrable Sturm involution $\sigma_g$ also corresponds to the full lap signature $\mathfrak{S}=\mathfrak{S}(T_g)$ of example \eqref{304}.

\sss{The frozen Vas tulip}\label{FT}

It remains to pass from the Neumann connection graph $\cCN$ of figure \ref{figNP}(left) to the delay connection graph of the Vas tulip $\CT$ in \ref{figNP}(right).
Vertical grading is by Morse index $i$; see \eqref{ij}, \eqref{iN}, \eqref{iP}.
For explanation of labels and symbols see figure \ref{figTNM}.
On the Neumann  left of figure \ref{figNP}, frozen equilibrium vertices are grouped in four min/max pairs $\{1,12\},\ \{3,11\},\ \{4,6\},\ \{8,10\}$, by their $q=4$ ``hidden'' frozen $\bO$ waves, alias ODE periodic orbits.
Each min/max pair corresponds to a 2-cycle of $\sigma_g$ in \eqref{sigmaT} with odd lap number $\ell=1$; see definition \ref{deflapsig} and the full lap signature \eqref{304}.

Passing from left to right, corollary \ref{corCPN}  determines the spatially reversible connection graph $\cC_g^\cP$ of PDE \eqref{101} under periodic boundary conditions.
Indeed, the relation $\sim$ merges the two vertices of each grouped frozen min/max pair into a single frozen wave.
Comparison with figure \ref{figtulip} shows that $\cC_g^\cP$ is indeed isomorphic to the connection graph $\CT$ of the Vas tulip $\AT$ of delay equation \eqref{501}, just as our heuristics have anticipated.
This completes our construction of the delay Vas tulip $\CT$ as a parabolic PDE connection graph $C_g^\cP$.

\bigskip \bigskip 

\noindent {\sc Acknowledgments.} 
Our very mode of mathematical thought owes much to the truly groundbreaking work of our dear jubilees and dedicatees.
Their friendly accessibility, their profound scholarship, and their books have shaped the thinking on dynamical systems of whole generations.
CR is deeply grateful for the continuous support by Isabel Rocha. 
Without her constant care, this endeavor would not have been possible. 
Mutual support and hospitality among all authors is also much cherished.
This work was partially supported by FCT/Portugal through the project 
UIDB/04459/2020 with DOI identifier 10-54499/UIDP/ 04459/2020, 
by the Taiwan National Science and Technology Council (NSTC) grant 113-2123-M-002-009,
by the Taiwan National Center of Theoretical Science (NCTS), and by the
Deutsche Forschungsgemeinschaft (DFG).

\end{document}